\documentclass[preprint,nopreprintline,10pt]{elsarticle}
\usepackage[margin=1.3in]{geometry}
\usepackage{amsmath,amssymb,amsthm,mathtools}
\usepackage{caption}
\usepackage{graphicx,subcaption}
\usepackage{tikz}
\usetikzlibrary{calc}
\usepackage{booktabs,multirow,makecell}
\usepackage{enumitem}
\setlist[enumerate]{label={(\arabic*)},leftmargin=*}
\usepackage{float}
\usepackage{hyperref}
\usepackage[nameinlink]{cleveref}
\crefname{appendix}{appendix}{appendices}
\Crefname{appendix}{Appendix}{Appendices}
\crefformat{appendix}{#2#1#3}
\Crefformat{appendix}{#2#1#3}

\allowdisplaybreaks

\let\b=\boldsymbol
\newcommand{\R}{\mathbb{R}}
\newcommand{\C}{\mathbb{C}}
\newcommand{\Z}{\mathbb{Z}}
\newcommand{\N}{\mathbb{N}}
\newcommand{\pa}{\partial}

\newenvironment{myproof}[1][\proofname]{%
  \begin{proof}[#1]\mbox{}\par\nobreak\ignorespaces
}{\end{proof}}

\theoremstyle{plain}
\newtheorem{theorem}{Theorem}[section]

\newtheorem{lemma}[theorem]{Lemma}

\theoremstyle{definition}
\newtheorem{definition}[theorem]{Definition}
\newtheorem{example}[theorem]{Example}
\newtheorem{remark}[theorem]{Remark}

\begin{document}

\begin{frontmatter}
\title{Quadrature by Matched Asymptotic Expansion (QBMAX): Accelerating Evaluation of Layer Potentials with Boundary Layers}
\begin{abstract}
Evaluating layer potentials with rapidly decaying kernels becomes
computationally challenging in the presence of boundary layers, as
occurs for Helmholtz layer potentials with large imaginary wavenumbers. We introduce
Quadrature by Matched Asymptotic Expansion (QBMAX), a high-order method that
extends the Quadrature by Expansion (QBX) framework by incorporating the
kernel's asymptotic decay behavior directly into the Taylor expansion of the
potential. Numerical results in both 2D and 3D show improvements of up to four digits over
QBX for the tested large-parameter cases, while the methods can perform similarly
for smaller parameters or more strongly curved geometries. A symbolic weighted
operation-count model indicates comparable costs for forming the two expansion
kernels. The Taylor expansions have the same formal order, while QBMAX reduces
the truncation-error constants in the flat-boundary model and generally exhibits
less convergence degradation in the reported experiments.
\end{abstract}

\author[inst1]{Chaoqi Lin}
\ead{chaoqi2@illinois.edu}
\author[inst1]{Xiaoyu Wei}
\ead{xywei@illinois.edu}
\author[inst1]{Andreas Kl\"ockner}
\ead{andreask@illinois.edu}
\address[inst1]{Siebel School of Computing and Data Science, University of Illinois at Urbana-Champaign}

\begin{keyword}
Modified Helmholtz equation \sep Integral equation method \sep Boundary layers
\sep Quadrature-by-expansion
\end{keyword}

\end{frontmatter}

\section{Introduction}
We consider the modified Helmholtz (or Yukawa) equation
\begin{equation}
    -\Delta u(\b x) + k^2 u(\b x) = 0, \qquad \b x \in \Omega, \label{eq:intro-modified-helmholtz}
\end{equation}
where $\Omega\subset\mathbb{R}^{n}$ is a bounded, adequately smooth domain with $n
\in \{2,3\}$, and $k>0$ is presumed to be large. Such equations arise in many
scientific and engineering applications, including thermoacoustic scattering
\cite{kirbyExactDomainTruncation2024}, the linearized Poisson-Boltzmann equation
\cite{kropinskiFastIntegralEquation2011}, and Navier-Stokes flows
\cite{afklintebergFastIntegralEquation2020}. For instance, in the semi-implicit
temporal discretization of the heat equation
\cite{quaifeFastIntegralEquation2011,kropinskiFastIntegralEquation2011c,fryklundIntegralEquationBased2020},
the parameter $k$ depends inversely on the time step. The method we present in
this paper generalizes to (unmodified) Helmholtz equations in the rapidly
evanescent regime, where $k \in \mathbb{C}$ with $\operatorname{Im}(k)\gg 1$. For
simplicity of the presentation, the bulk of this article is written to focus on
\eqref{eq:intro-modified-helmholtz}. In \Cref{ex:helmholtz-starfish}, we present a
numerical experiment showing the applicability of the method to the Helmholtz
equation.

We focus on integral equation approaches to solve~\eqref{eq:intro-modified-helmholtz}.
Integral equation methods make use of the fundamental solution of the PDE. In the
Yukawa case, the Green's function is given by
\[
    \mathcal{G}_k(\b x,\b y) = \begin{cases}
         \dfrac{1}{2\pi} K_0(k|\b x-\b y|), & n=2, \\
         \dfrac{1}{4\pi} \dfrac{e^{-k|\b x-\b y|}}{|\b x-\b y|}, & n=3,
    \end{cases}
\]
where $K_0$ is the modified Bessel function of the second kind of order zero.
Depending on the boundary conditions, the solution $u$ can be represented as a
linear combination of layer potentials:
\begin{equation}
 u (\b x) = \alpha \int_{\partial \Omega} \mathcal{G}_k(\b x,\b y)\,\sigma(\b y)\,dS_{\b y} + \beta \int_{\partial \Omega} \frac{\partial \mathcal{G}_k(\b x,\b y)}{\partial \b \nu_{\b y}}\,\sigma(\b y)\,dS_{\b y},\quad \b x \in \mathbb{R}^n \setminus \pa \Omega, \label{eq:layer-potential-form}
\end{equation}
for some $\alpha, \beta \in\mathbb{R}$ and an unknown density $\sigma$. A direct
advantage of such integral formulations is the dimension reduction: only
boundary unknowns are required. Moreover, the density $\sigma$ is often
determined via a second-kind Fredholm equation that typically has benign conditioning (see,
e.g., \cite{kropinskiFastIntegralEquation2011} and
\cite{jiangSecondKindIntegral2013c}). However, the kernel $\mathcal{G}_k(\b x,\b y)$
is singular when $\b x=\b y$, and the evaluation of the layer potential requires
specialized treatment when the target point $\b x$ is close to the source geometry $\partial \Omega$ and often
requires special numerical schemes.
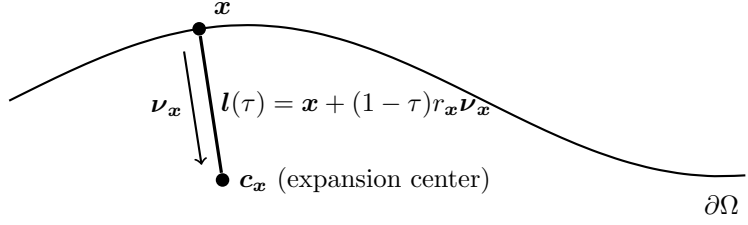
\begin{figure}
    \centering
    \begin{tikzpicture}
        \draw[domain=0:3.1*pi, smooth, variable=\x, black, thick]
             plot ({\x}, {sin(\x / 2 r)});

        \node (A) at ({3.2*pi/4}, {sin(3.2*pi/8 r)}) [circle, fill, inner sep=1.8pt, label=above right:{$\b x$}] {};
        \node (B) at ({3.2*pi/4 + 2*1/2 * cos(3.2*pi/8 r)}, {sin(3.2*pi/8 r)-2})
        [circle, fill, inner sep=1.8pt, label=right:{$\b c_{\b x} \ (\text{expansion center})$}] {};
        \draw[-, very thick, black] (A) -- node[midway, right] {$\b l(\tau) =\b x  +  (1-\tau) r_{\b x} \b \nu_{\b x}$} (B);

\coordinate (A_shifted) at ({3.2*pi/4 - 0.2}, {sin(3.2*pi/8 r) - 0.3});
\coordinate (B_shifted) at ({3.2*pi/4 - 0.2 +  1.5*1/2 * cos(3.2*pi/8 r)}, {sin(3.2*pi/8 r) - 0.3 - 1.5});
        \draw[->, thick, black] (A_shifted) -- (B_shifted) node[midway, left] {$\b \nu_{\b x}$};
        \node[below] at ({3*pi}, {sin((3*pi)/2 r)-1/8}) {$\partial \Omega$};
    \end{tikzpicture}
    \caption{QBX expansion near curve $\partial \Omega$ at target $\b x$, with
    expansion center $\b c_{\b x} = \b x + r_{\b x} \b \nu_{\b x}$ at distance $r_{\b x}$, where $\b \nu_{\b x}$
    is the outward unit normal at $\b x$. The expansion line $\b l(\tau) = \b x +
    (1-\tau) r_{\b x} \b \nu_{\b x}$ connects the target to the expansion center.}
    \label{fig:qbx_line_expansion}
\end{figure}

In this work, we focus on the near-field evaluation (i.e., when the target point $\b x$ is close to
or on the boundary). For the set of problems under consideration, that is the
dominant concern, since the layer potential rapidly decays to zero elsewhere. To
that end, we consider a variant of the Quadrature by Expansion (QBX) method
\cite{klocknerQuadratureExpansionNew2013, walaFastAlgorithmQuadrature2019,
afklintebergAdaptiveQuadratureExpansion2018}. A marked advantage of QBX is its
ability to handle a broad cross-section of kernels, dimensionalities, and
singularities within one unified, straightforward framework, while also
unifying singular (on-surface) and near-singular (off-surface) evaluation.

By forming off-surface local expansions of the layer potentials, QBX effectively
`smooths out' the kernel singularities, allowing for accurate near-field
evaluations via conventional high-order quadrature rules. Suppose $\b x \in
\partial \Omega$ with an expansion line $\b l(\tau)$ as defined in
\Cref{fig:qbx_line_expansion}. The value $u$ along the expansion line $\b l$ can be
approximated by a $p^{\text{th}}$-order Taylor expansion centered at $\tau = 0$,
\begin{equation}
\begin{split}
(u \circ \b l) (\tau)  &\approx \mathcal{T}_{p}[u \circ \b l ](\tau):=\sum_{m=0}^{p} \frac{\partial_{\tau}^{m} (u \circ \b l)(\tau) \big|_{\tau=0}}{m!} \tau^{m} \\
&=  \int_{\partial \Omega} \underbrace{\sum_{m=0}^{p} \frac{\partial_{\tau}^m [\tau  \mapsto  \mathcal{K}(\b l(\tau), \b y)]\big|_{\tau=0}}{m!} \tau ^{m} }_{\text{expansion kernel}} \sigma(\b y)\,dS_{\b y},
\end{split}
\label{eq:qbx-line-taylor-intro}
\end{equation}
where $\mathcal{K}$ denotes a kernel function such as $\mathcal G_k$ and $\mathcal{T}_{p}$ denotes a
$p+1$ term Taylor expansion about $\tau = 0$. However, for rapidly-decaying Yukawa potentials, the QBX line expansion
\eqref{eq:qbx-line-taylor-intro} can have a large truncation error relative to the
magnitude of the potential. This occurs for large parameter values $k$, where
$\mathcal{K}$ decays exponentially in $k|\b x-\b y|$ and a fixed-order Taylor
polynomial need not capture that decay over a fixed expansion radius. Analogously, the Taylor
expansion of $e^{z}$ about $z=0$,
\[
  e^{z} = \sum_{j=0}^{p} \frac{z^{j}}{j!} + R_{p}(z)
\]
has remainder $|R_{p}(z)| \leq e^{z}\,|z|^{p+1}/(p+1)!$, which grows with
$z \geq 0$ for fixed expansion order $p$.

To mitigate the possibly large truncation error, we form a new series expansion by encoding the asymptotic decay
information to achieve a better approximation. 
We do so in the context of the QBX framework by introducing a smooth
asymptotic scaling function
\[u_{\text{ma}}: [0,1] \to \mathbb{R}, \quad u_{\text{ma}}(\tau) = \exp(k|\b l(\tau)-\b x|)\]
that removes the leading exponential decay of the potential along the expansion
line before Taylor approximation.
Then, $(u \circ \b l) (\tau)$ can be approximated by
\begin{equation}
(u \circ \b l) (\tau) \approx \frac{\mathcal{T}_{p}[(u \circ \b l)  \cdot u_{\textup{ma}} ] (\tau)}{u_{\textup{ma}}(\tau)}.
\label{eq:qbx-line-taylor-modified}
\end{equation}
The Taylor remainder theorem guarantees that the asymptotic behavior of the
truncation error necessarily remains the same as that of \eqref{eq:qbx-line-taylor-intro}.
Since $u_{\text{ma}}$ is specific to a given target $\b x$, the expansion in~\eqref{eq:qbx-line-taylor-modified}
is also  \emph{target-specific} in the sense of~\cite{siegelLocalTargetSpecific2018}.
We term this approach
\emph{Quadrature by Matched Asymptotic Expansion} (QBMAX). 

Related ideas making use of modified Taylor
expansions have been used in other contexts, e.g., singularity removal
\cite{lamonRemovalSingularitiesTaylor1989} and function approximations
\cite{bulosModifiedTaylorsApproximation2005,asschePadeHermitePadeApproximation2006a,howardDualTaylorSeries2019},
but, to the best of our knowledge, their application in the contexts
of rapidly decaying kernels and integral equations is novel.

It is natural to wonder about the error behavior of \eqref{eq:qbx-line-taylor-modified}
in relation to the original expansion.
Since the asymptotic behavior is the same, any error results must account
for constant factors. We provide a detailed analysis of the constant factors in the truncation error
of \eqref{eq:qbx-line-taylor-modified} via a salient example in \Cref{sec:flat-panel-truncation-error}.
There, we also study the factor by which the error is decreased compared to the
original QBX expansion. We furthermore present numerical evidence
for less idealized settings in \Cref{sec:numerical-results}.
The infinite-flat-boundary model in \Cref{sec:flat-panel-truncation-error} and
the experiments in \Cref{sec:numerical-results} indicate that the
transformation~\eqref{eq:qbx-line-taylor-modified} can substantially reduce the
relative truncation error for large $k$.

The modified Helmholtz equation with large $k$ is a well-studied
problem. Ockendon and Tew \cite{ockendonThinLayerSolutionsHelmholtz2012}
provided a detailed overview of the ``thin-layer'' structures of the solutions
in two dimensions. For close evaluation of the layer potentials, the numerical
difficulty is twofold: the kernel is singular or nearly singular near the
boundary, and the kernel varies rapidly in the vicinity of the singularity.
Fryklund et al.\ \cite{fryklundAdaptiveKernelsplitQuadrature2022} extend the
kernel-split method \cite{helsingEvaluationLayerPotentials2008a} to a class of
parameter-dependent modified PDEs, including the Yukawa equation. Their method
addresses the rapid decay of the kernel by recursively refining panels near the
target point until the kernel is sufficiently resolved. It applies to a wide
range of $k$ values, with additional computational cost scaling as
$\mathcal{O}(\log k)$ in two dimensions.

Recently, Stein and Barnett \cite{steinQuadratureFundamentalSolutions2021} introduced
Quadrature by Fundamental Solutions (QFS), a kernel-independent method that represents
a given layer potential by an equivalent off-surface source distribution.
Although QFS is related in spirit to the classical method of fundamental solutions
\cite{linFastSolutionThreeDimensional2016,eiMethodFundamentalSolutions2022},
the effective sources are constructed to approximate the original potential,
rather than to solve a boundary value problem directly. In the subsequent
work \cite{steinSpectrallyAccurateSolutions2022}, QFS-type close evaluation
was incorporated into a solver for inhomogeneous modified Helmholtz equations
at moderately large $k$. There, the original boundary integral operators are
discretized by a Nystr\"om method with high-order Alpert quadrature
\cite{alpertHybridGaussTrapezoidalQuadrature1999}, and the source distribution of QFS
is determined by matching the potential at boundary collocation points. QFS
avoids singular quadrature rules and reduces the cost in the subsequent
evaluations, but the Nystr\"om discretization must resolve
the kernel properly (for example, the panel size near the singularity should be $\mathcal{O}(1/k)$)
in the first place.

Other related approaches include density interpolation methods
\cite{fariaGeneralpurposeKernelRegularization2021}, special quadrature rules
\cite{kapurHighOrderCorrectedTrapezoidal1997,wuZetaCorrectionNew2021}, and modified
layer-potential representations \cite{carvalhoModifiedRepresentationsClose2021a}.
These methods primarily address the kernel singularity; handling the rapid
kernel variation would require combining them with local adaptive
refinement.

QBMAX includes the leading asymptotic decay in the line expansion. In the
reported experiments this permits the same refinement sequences to be used for
all tested values of $k$; we do not claim that the required resolution is
uniform in $k$ for arbitrary geometries or parameter ranges.
It is appealing because it provides a framework for high-order
numerical methods that straightforwardly generalizes to (at least) two and three
dimensions as well as to a variety of kernels, including both oscillatory and
non-oscillatory ones. The method focuses particularly on near- and on-surface
evaluations in the setting of ``thin-layer'' solutions. The range of PDEs to which QBMAX extends
naturally includes the modified biharmonic equation
\cite{jiangSecondKindIntegral2013c}, the modified Stokes equations
\cite{afklintebergFastIntegralEquation2020}, and the Yukawa--Beltrami equations
\cite{kropinskiIntegralEquationMethods2016}.

In this work, we make the following contributions:
\begin{enumerate}
    \item We introduce QBMAX, a high-order method that accelerates the conventional
          QBX method by including the kernel's exponential decay behavior
          directly into the Taylor expansion. This modification significantly
          improves accuracy for many of the tested large values of $k$. A
          symbolic weighted operation-count model indicates a cost for forming
          the QBMAX expansion kernel comparable to that for QBX. Across each
          reported refinement study we use the same meshes for all values of
          $k$; within these tests, QBMAX is substantially less sensitive to $k$
          than QBX in the rapidly decaying regime.

    \item In \Cref{sec:flat-panel-truncation-error}, we provide a detailed
          truncation error analysis for an infinite flat boundary in 2D. The analysis
          explicitly shows the reduction in constant factors in the error
          estimates and facilitates understanding of error behavior in more
          general settings.

    \item In the two- and three-dimensional numerical examples presented,
          QBMAX performs similarly to QBX in some small-parameter or curved-
          geometry cases and improves accuracy by up to four digits in some
          large-parameter cases. Its observed convergence rates generally
          degrade less with increasing $k$, although the rates vary with the
          geometry.
\end{enumerate}

\section{Fundamental Solutions and Layer Potentials}
\subsection{Fundamental Solutions}
\begin{figure}[H]
	\centering \resizebox{0.7\textwidth}{!}{%
		\input{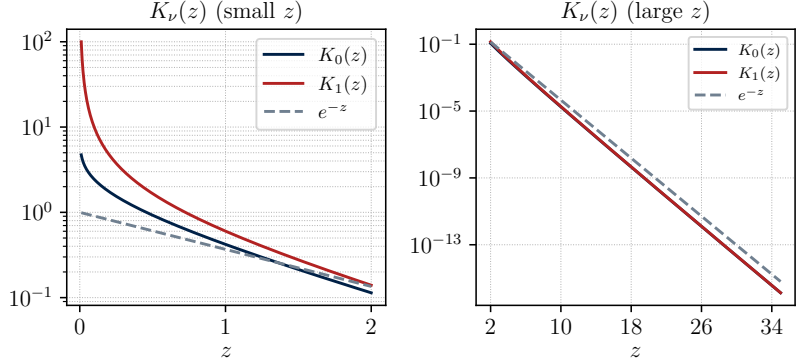}
	} \caption{Plots of the modified Bessel functions $K_0(z)$ and $K_1(z)$.}
	\label{fig:bessel}
\end{figure}

For $k>0$ and $n \in \{2,3\}$, the fundamental solutions $\mathcal{G}_k(\b x,\b y)$
to the modified Helmholtz equation
 \[
	(-\Delta_{\b x} + k^2)\mathcal{G}_k(\b x,\b y) = \delta(\b x-\b y), \quad \b x,\b y \in \mathbb{R}^{n}
,\]
are given by
\[
	\mathcal{G}_k(\b x,\b y) = \begin{cases}
		 \frac{1}{2\pi} K_0(k|\b x-\b y|), & n=2, \\
	\frac{1}{4\pi}	\frac{e^{-k|\b x-\b y|}}{|\b x-\b y|}, & n=3,
	\end{cases}
\]
where $K_0(k|\b x-\b y|)$ is the modified Bessel function of the second kind of order
zero. For $z>0$, $K_0(z)$ admits the expansion
\cite[\href{https://dlmf.nist.gov/10.31.E2}{(10.31.2)}]{DLMFNISTDigital},
\begin{equation*}
	\begin{split}
	K_0(z) = \ln\!\bigl(\tfrac{2}{z}\bigr)
		&+ \ln\!\bigl(\tfrac{2}{z}\bigr)\sum_{n=1}^{\infty}\frac{1}{(n!)^2}\Bigl(\frac{z^2}{4}\Bigr)^{n}
		- \gamma\Bigl(1 + \sum_{n=1}^{\infty}\frac{1}{(n!)^2}\Bigl(\frac{z^2}{4}\Bigr)^{n}\Bigr) \\
		&+ \sum_{n=1}^{\infty}\!\Bigl(1+\tfrac12+\cdots+\tfrac1n\Bigr)\frac{1}{(n!)^2}\Bigl(\frac{z^2}{4}\Bigr)^{n},
	\end{split}
\end{equation*}
where $\gamma$ is the Euler-Mascheroni constant. The leading singular behavior
of $K_0(z)$ and $\partial_z K_0(z)$ near $z=0$ is $-\ln(z)$ and $-1/z$,
respectively. These singularities have been studied extensively in the QBX
setting~\cite{klocknerQuadratureExpansionNew2013,
epsteinConvergenceLocalExpansions2013a,klintebergErrorEstimationQuadrature2017}.
However, $K_\nu(z)\sim\sqrt{\frac{\pi}{2z}}\,e^{-z}$ as $z\to\infty$ with fixed
order $\nu$ \cite[\href{https://dlmf.nist.gov/10.40.E2}{(10.40.2)}]{DLMFNISTDigital}. Thus,
$K_\nu(k|\b x-\b y|)$ decays rapidly away from $\b x=\b y$ and, for large $k$, is
concentrated in a thin region around the singularity (see~\Cref{fig:bessel}).

\subsection{Layer Potentials}\label{sec:layer-potentials}
Let $\Omega \subset \R^{n}$ be a bounded domain with smooth
boundary $\pa \Omega$. We define the single and double layer potentials as
follows.
\begin{definition}[Single layer potential]
Assume $\sigma \in L^2(\pa \Omega)$. The single layer potential is defined as
\begin{equation}
	\mathcal{S}_{k}\sigma(\b x) = \int_{\pa \Omega} \mathcal{G}_k(\b x,\b y) \sigma(\b y) dS_{\b y}, \quad \b x \in \mathbb{R}^{n}\setminus\partial\Omega \label{eq:single-layer}
\end{equation}
\end{definition}
\begin{definition}[Double layer potential]
	Assume $\sigma \in L^2(\pa \Omega)$. The double layer potential is defined as
	\begin{equation}
		\mathcal{D}_{k}\sigma(\b x)  = \int_{\pa \Omega} \frac{\partial \mathcal{G}_k(\b x,\b y)}{\partial \b \nu_{\b y}} \sigma(\b y) dS_{\b y}, \quad \b x \in \mathbb{R}^{n}\setminus\partial\Omega, \label{eq:double-layer}
	\end{equation}
	where $\b \nu_{\b y}$ is the unit outward normal to $\pa \Omega$ at $\b y$ and $\frac{\pa
	\mathcal{G}_{k}(\b x,\b y)}{\pa \b \nu_{\b y}}:= \langle \nabla_{\b y} \mathcal{G}_{k}(\b x,\b y),
	\b \nu_{\b y}\rangle.$
\end{definition}
Often, $\b x$ is referred to as the \textit{target} point, representing the
location where the potential is being evaluated, and $\b y$ as the \textit{source}
point, representing the location on the boundary $\partial \Omega$ contributing
to the potential. As $\mathcal{G}_k$ has the same leading singularity as the
Laplace kernel, $\mathcal{S}_{k}\sigma$ is continuous on $\R^{n}$ and the jump
relations for $\mathcal{D}_{k}\sigma$ and $\mathcal{S}'_{k}\sigma$ (e.g.,
\cite{quaifeFastIntegralEquation2011},
\cite{hsiaoBoundaryIntegralEquations2021}) are as follows
\begin{align*}
	\mathcal{D}_{k} \sigma (\b x) &= \lim_{h \searrow 0^{+}} \mathcal{D}_{k}  \sigma (\b x \mp h \b \nu_{\b x}) \pm \frac{1}{2} \sigma(\b x) , \; \text{ for } \b x\in \pa \Omega, \\
	\mathcal{S}'_{k} \sigma (\b x) &= \lim_{h \searrow 0^{+}} \mathcal{S}'_{k}  \sigma (\b x \mp h \b \nu_{\b x}) \mp \frac{1}{2} \sigma(\b x) , \; \text{ for } \b x\in \pa \Omega,
\end{align*}
where $\mathcal{S}'_{k}\sigma(\b x\mp h \b \nu_{\b x}):= \frac{\pa}{\pa
\b \nu_{\b x}}\int_{\partial \Omega}\mathcal{G}_{k}(\b x\mp h \b \nu_{\b x}, \b y) \sigma(\b y)dS_{\b y}$.

\subsection{Quadrature by Expansion and Taylor Expansions of Layer Potential on a Line}\label{sec:qbx-line-taylor}
The QBX method \cite{klocknerQuadratureExpansionNew2013} provides a systematic
and generalizable approach for the close evaluation of layer potentials with
singular kernels. When the target point is away from the boundary, layer
potentials are smooth and can be integrated accurately using standard quadrature
rules (e.g., Gauss-Legendre, trapezoidal rule). However, as the target
approaches the boundary, standard quadrature rules suffer from rapidly growing
errors due to kernel singularities.

QBX can be viewed as a kernel regularization method that uses Taylor expansions
centered at points ``away'' from the boundary to approximate targets in the near
field. Then, the QBX expansion kernels \eqref{eq:qbx-line-taylor-intro} can be
integrated accurately with standard quadrature rules. This approach relies on
the analyticity of layer potentials up to the boundary under appropriate
regularity assumptions (see~\Cref{lem:regularity-of-layer-potentials}). We refer
interested readers to
\cite{epsteinConvergenceLocalExpansions2013a,coltonInverseAcousticElectromagnetic2013,dautrayMathematicalAnalysisNumerical1999}
for additional details on the regularity of layer potentials.
\begin{lemma}[Regularity of layer potentials~\cite{dautrayMathematicalAnalysisNumerical1999}] \label{lem:regularity-of-layer-potentials}
   Suppose $\Omega$ is of class $\mathcal{C}^{m+1,\beta}$ and $\sigma$ is of
   class $\mathcal{C}^{m,\beta}$ for some $m \in \mathbb{N}$ and $0 < \beta <
   1$. Then:
   \begin{enumerate}
       \item $\mathcal{S}_{k}\sigma$ is of class $\mathcal{C}^{m+1,\beta}$ on
             $\overline{\Omega}$ and $\mathbb{R}^n \setminus \Omega$.
       \item $\mathcal{D}_{k}\sigma$ is of class $\mathcal{C}^{m,\beta}$ on
             $\overline{\Omega}$ and $\mathbb{R}^n \setminus \Omega$.
   \end{enumerate}
\end{lemma}
Next, we summarize the line QBX expansion, a variant of the original QBX
method~\cite{klocknerQuadratureExpansionNew2013}, in the following steps. Fix a
target point $\b x \in \partial \Omega$ and define:
\begin{enumerate}
    \item \textbf{Expansion centers:} Define an off-surface expansion center
          $\b c_{\b x}$ on the normal line at $\b x$ with radius $r_{\b x}:=|\b x-\b c_{\b x}|$.

    \item \textbf{Expansion line:} The expansion line $\b l(\tau)$ is defined as
          $\b l(\tau) = \b c_{\b x} + \tau (\b x - \b c_{\b x})$ for $\tau \in [0,1]$ (see
          \Cref{fig:qbx_line_expansion}).

    \item \textbf{Line Taylor expansion:} Recall from \eqref{eq:qbx-line-taylor-intro} that $(u \circ \b l)(\tau)$ can be approximated by a
          $p^{\text{th}}$-order Taylor expansion centered at $\tau = 0$:
    \begin{equation}
        (u \circ \b l)(\tau) \approx \mathcal{T}_{p}[u \circ \b l](\tau) = \int_{\partial \Omega} \mathcal{K}_{\text{QBX}}^{p}(\b l(\tau), \b y)\,\sigma(\b y) \, dS_{\b y}, \label{eq:line-qbx-taylor-general}
    \end{equation}
    where the expansion kernel $\mathcal{K}_{\text{QBX}}^{p}$ is defined as
    \begin{equation}
        \mathcal{K}_{\text{QBX}}^{p}(\b l(\tau), \b y) := \sum_{m=0}^{p} \frac{\partial_{\tau}^{m} \left[\tau \mapsto \mathcal{K}(\b l(\tau), \b y)\right] \big|_{\tau = 0}}{m!} \tau^{m}, \label{eq:regularized-qbx-kernel}
    \end{equation}
    and $\mathcal{K}$ is the kernel function of $u$.
\end{enumerate}
The truncation error follows from Taylor's remainder theorem
\begin{equation}
	|(u \circ \b l)(\tau) - \mathcal{T}_p[u \circ \b l](\tau)|
	\leq \frac{(r_{\b x}\tau)^{p+1}}{(p+1)!}\sup_{s \in [0,\tau]} \left|\big(\b n(\b x)\cdot \nabla_{\b z}\big)^{p+1} u(\b z)\big\rvert_{\b z = \b l(s)}\right|, \label{eq:QBX-truncation}
\end{equation}
where $\b n(\b x):=(\b x - \b c_{\b x})/|\b x - \b c_{\b x}|$.
There are two types of errors in the QBX line method:
truncation error and quadrature error. The latter arises from numerical
integration of the regularized kernel \eqref{eq:line-qbx-taylor-general} and
has been studied in
\cite{klocknerQuadratureExpansionNew2013,
epsteinConvergenceLocalExpansions2013a,klintebergErrorEstimationQuadrature2017,siegelLocalTargetSpecific2018}.
As shown in
\cite{klocknerQuadratureExpansionNew2013,epsteinConvergenceLocalExpansions2013a},
the QBX quadrature error can be estimated as follows:

\begin{theorem}[QBX quadrature error, \cite{klocknerQuadratureExpansionNew2013, epsteinConvergenceLocalExpansions2013a}]\label{thm:qbx-series-quadrature-error}
    Suppose $\partial \Omega \subset \mathbb{R}^{2}$ is a smooth closed curve
    and $B_{r}(\b c)$ is an open ball centered at $\b c$ with radius $r$ such that
    $\overline{B_{r}(\b c)} \cap \partial \Omega = \{ \b x \}$. Further, suppose the
    boundary $\partial \Omega$ is divided into panels of uniform size $h$, and a
    $q$-point composite Gauss--Legendre quadrature rule is used to integrate the
    $p^{\text{th}}$-order QBX expansion kernels
    \cite[eq.~9]{klocknerQuadratureExpansionNew2013}. For $0 < \beta < 1$, there
    exists a constant $C_{p, q, \beta}$ such that if
    $\sigma \in \mathcal{C}^{2q,\beta}(\partial\Omega)$, then the quadrature
    error of the QBX method satisfies
    \begin{equation}
       |\text{Quadrature Error}| \leq C_{p,q,\beta} \left(\frac{h}{4r}\right)^{2q} \left\|\sigma\right\|_{\mathcal{C}^{2q, \beta}(\partial \Omega)}.    \label{eq:quadrature-error-estimate}
    \end{equation}
\end{theorem}

The above estimates hold for the single layer potentials, with extensions to
three dimensions available
\cite{klintebergErrorEstimationQuadrature2017,siegelLocalTargetSpecific2018}. In
the original QBX method \cite{klocknerQuadratureExpansionNew2013}, the expansion
kernels are defined in terms of volumetric series expansions rather than line expansions.
For completeness, we also provide the quadrature error estimates for the line
QBX expansion in \Cref{sec:appendix-quad-error}.

The expansion radius $r$ and expansion order $p$ control the truncation error.
However, choosing $r$ and $p$ is not straightforward. Decreasing $r$ or
increasing $p$ requires higher quadrature resolution as the kernel becomes more
singular. Based on \eqref{eq:quadrature-error-estimate}, $r$ should be chosen so
that $\frac{h}{4r} < 1$ to control the quadrature error up to a finite
tolerance. In practice, the restriction $\frac{h}{4r} < 1$ can be relaxed, and
finding a good balance between these parameters is crucial for achieving both
accuracy and efficiency \cite{klocknerQuadratureExpansionNew2013,
rachhFastAlgorithmsQuadrature2017, walaFastAlgorithmQuadrature2019}.

In a concrete realization of QBX, care must be taken to ensure
a number of further properties. We provide a brief overview of these
following~\cite{walaFastAlgorithmQuadrature2019}. First, the source geometry is
refined so that the expansion disk at each target point does not interfere with
the source geometry. The result of this is referred to as \textit{stage-1 discretization}.
Next, adaptive refinement is performed to ensure sufficient quadrature
resolution for each target point. The result of this is termed \textit{stage-2
discretization}. If needed, \textit{stage-2-quad discretization} can be
generated by upsampling the quadrature in \textit{stage-2} to help resolve the
kernel. Once we have adequate source discretization, \eqref{eq:line-qbx-taylor-general}
can be evaluated with controlled \textit{asymptotic error}.

The kernel-independent nature of the QBX method can be two-edged: it is
generally applicable to a broad class of singular kernels, but the asymptotic
error estimates do not make the constants explicit. As the numerical examples
below illustrate, the QBX relative error and its pre-asymptotic convergence behavior can
deteriorate as $k$ increases, making the method substantially less efficient at
achieving a prescribed relative accuracy in the rapidly decaying regime. In the following
sections, we will discuss how to reduce the truncation error constants and
improve the convergence of QBX for large $k$.

\section{Quadrature by Matched Asymptotic Expansion}\label{sec:qbmax}
Our goal in this section is to provide an understanding of the constant factors
in the truncation error estimates. Unfortunately, obtaining estimates with
explicit constants for general geometries poses significant difficulty,
particularly due to chain rule terms involving high-order derivatives of the
geometry. Consequently, we present a careful analysis of the truncation error
behavior for an infinite flat boundary in 2D. For large values of $k$, the kernel
$\mathcal{G}_k$ decays rapidly away from the singularity, making the tail
contribution from the integral negligible. Therefore, as an illustrative example,
we consider a source domain of $\mathbb{R}\times \{0\}$, with
Fourier-mode densities.

\subsection{Infinite Flat Boundary: Truncation Error Analysis}\label{sec:flat-panel-truncation-error}
Assume $\b x=(x_1,x_2), \b y = (y_1,y_2) \in \mathbb{R}^2$ and $\sigma(\b y) = e^{-i \omega y_1}$ for some $\omega \in \mathbb{R}$. Observe
\begin{align*}
\int_{\mathbb{R} \times \{0\}} K_0(k|\b x-\b y|) \sigma(\b y) \,dS_{\b y} &= \int_{\mathbb{R}} K_0(k\sqrt{(x_1 - y_1)^2 + x_2^2}) e^{-i \omega y_1} \,dy_1 \\
&= \int_{\mathbb{R}} K_0(k\sqrt{s^2 + x_2^2}) e^{-i \omega (x_1 - s)} \,ds \\
&= e^{-i \omega x_1} \int_{\mathbb{R}} K_0(k \sqrt{s^2 + x_2^2}) e^{i \omega s} \,ds.
\end{align*}
By \Cref{lem:flat-panel}, one shows
\begin{equation}
\int_{\mathbb{R} \times \{0\}} K_0(k|\b x-\b y|) \sigma(\b y) \,dS_{\b y} = e^{-i \omega x_1} \frac{\pi}{\sqrt{k^2 + \omega^2}} \exp(-|x_2| \sqrt{k^2 + \omega^2}).
\label{eq:infinite-flat-boundary-solution}
\end{equation}
Below, we focus on the error analysis for the on-surface evaluation as the
truncation error of \eqref{eq:line-qbx-taylor-general} is (typically) the largest when
$\tau=1$. Assume $\b x^{*} = (x_1,0)$ with expansion center $\b c = (x_1,r)$ for some
$r > 0$. Denote $\b l(\tau)= \b c + \tau (\b x^{*}-\b c)$. The $p^{\text{th}}$-order QBX
line expansion of $(u \circ \b l)(\tau)$ centered at $\tau = 0$ and evaluated at
$\tau =1$ is given by
\begin{equation}
\begin{split}
 \mathcal{T}_{p}[u \circ \b l](1) &= e^{-i \omega x_1} \frac{\pi}{\sqrt{k^2 + \omega^2}} \sum_{m=0}^{p} \frac{\partial_\tau^m[\exp(- r(1-\tau) \sqrt{k^2 + \omega^2})]}{m!} \bigg\rvert_{\tau=0} \\
&= e^{-i \omega x_1} \frac{\pi}{\sqrt{k^2 + \omega^2}} \exp(-r\sqrt{k^2 + \omega^2}) \sum_{m=0}^{p} \frac{(r\sqrt{k^2 + \omega^2})^m}{m!}
\end{split}
\label{eq:qbx-flat-panel-expansion}
\end{equation}
with truncation error
\begin{equation}
|u(\b x^{*}) - \mathcal{T}_{p}[u \circ \b l](1)| =  \frac{\pi}{\sqrt{k^{2} + \omega^{2}}} \frac{\gamma(p+1, r\sqrt{k^{2}+ \omega^{2}} )}{\Gamma(p+1, 0)},
\label{eq:qbx-flat-panel-exact-truncation-error}
\end{equation}
where $\gamma$ denotes the lower incomplete gamma function and $\Gamma$ the gamma function, so that $\Gamma(p+1,0) = \Gamma(p+1) = p!$.

\begin{figure}[H]
    \centering \resizebox{0.5\textwidth}{!}{%
	\input{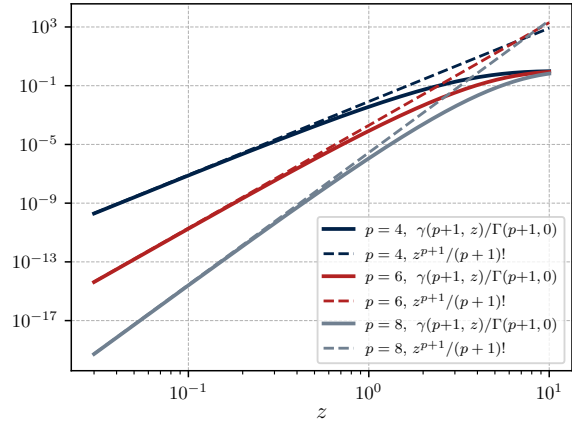}
} \caption{Plots of $\frac{\gamma(p+1,z)}{\Gamma(p+1,0)}$ and
$\frac{z^{p+1}}{(p+1)!}$.} \label{fig:gamma_bound}
\end{figure}

Additionally, the Taylor remainder theorem shows that the
truncation error can also be bounded as
\begin{equation}
 	|u(\b x^{*}) - \mathcal{T}_{p}[u \circ \b l](1)| \leq  \frac{\pi}{\sqrt{k^{2} + \omega^{2}}} \frac{(r\sqrt{k^{2} + \omega^{2}})^{p+1}}{(p+1)!} \label{eq:QBX-flat-panel-truncation-error}
\end{equation}
After division by the on-surface magnitude $|u(\b x^*)|=
\pi/\sqrt{k^2+\omega^2}$, the exact relative error in
\eqref{eq:qbx-flat-panel-exact-truncation-error} is the regularized incomplete
gamma function shown in \Cref{fig:gamma_bound}. For fixed $p$, this relative
error approaches one when $r\sqrt{k^2+\omega^2}$ becomes large, even though the
absolute error is asymptotic to $\pi/k$ for fixed $r$ and $\omega$ as
$k\to\infty$. Maintaining a prescribed relative truncation accuracy therefore
requires $r\sqrt{k^2+\omega^2}$ to remain small (typically below one in this
model). Decreasing $r$ achieves this but requires higher quadrature resolution,
leading to increased computational cost for large $k$, particularly in 3D.
One possible resolution is to form a line expansion that includes the asymptotic
information of the layer potential. The exact solution \eqref{eq:infinite-flat-boundary-solution} decays in
$|x_2|$ at rate $\sqrt{k^2+\omega^2}$. QBMAX removes the leading large-$k$ part
$k$ of this rate, leaving the residual rate $\sqrt{k^2+\omega^2}-k$; this is
zero for a constant density ($\omega=0$) and small for fixed $\omega$ as
$k\to\infty$. Accordingly, instead of forming a vanilla line expansion of $(u
\circ \b l) (\tau)$ at $\tau= 0$, we form the Taylor expansion of 
\[
(u \circ \b l) (\tau) \cdot \exp(k|\b l(\tau) - \b x^{*}|)\quad \text{at } \tau = 0
\]
and evaluate the resulting expansion at $\tau=1$.
To recover the original layer potential, the new expansion needs to be scaled by
$\exp(-k|\b l(1) - \b x^{*}|)=1$. Precisely, $u(\b x^{*})$ is approximated by
\begin{equation}
\begin{split}
u(\b x^{*}) & \approx \exp(-k|\b l(1) - \b x^{*}|) \sum_{m=0}^{p} \frac{\pa_\tau^{m} \left[(u \circ \b l)(\tau) \cdot \exp(k |\b l(\tau) - \b x^{*}|) \right] }{m!} \bigg\rvert_{\tau=0} \\
&= e^{-i \omega x_1} \frac{\pi}{\sqrt{k^2 + \omega^2}} \sum_{m=0}^{p} \frac{\partial_\tau^m[\exp(- r(1-\tau) \sqrt{k^2 + \omega^2}) \exp(r(1-\tau)k)]}{m!}  \bigg\rvert_{\tau=0}  \\
&= e^{-i \omega x_1} \frac{\pi}{\sqrt{k^2 + \omega^2}} \exp(-r(\sqrt{k^2 + \omega^2} - k)) \sum_{m=0}^{p} \frac{\partial_\tau^m[\exp(\tau r(\sqrt{k^2 + \omega^2} -k))]}{m!} \bigg\rvert_{\tau=0}  \\
&= e^{-i \omega x_1} \frac{\pi}{\sqrt{k^2 + \omega^2}} \exp(-r(\sqrt{k^2 + \omega^2} - k)) \sum_{m=0}^{p} \frac{(r(\sqrt{k^2 + \omega^2} -k))^m}{m!}
\end{split}
\label{eq:qbmax-flat-panel-expansion}
\end{equation}
with truncation error
\begin{equation}
\frac{\pi}{\sqrt{k^{2} + \omega^{2}}} \frac{\gamma(p+1, r(\sqrt{k^{2}+ \omega^{2}} - k) )}{\Gamma(p+1, 0)} \label{eq:qbmax-flat-panel-truncation-error}.
\end{equation}
Remarkably, when $\omega = 0$, the new expansion \eqref{eq:qbmax-flat-panel-expansion} has truncation error $0$. That is, for a constant layer density, it is
\emph{exact}.

\begin{figure}[H]
    \centering \resizebox{\textwidth}{!}{%
        \input{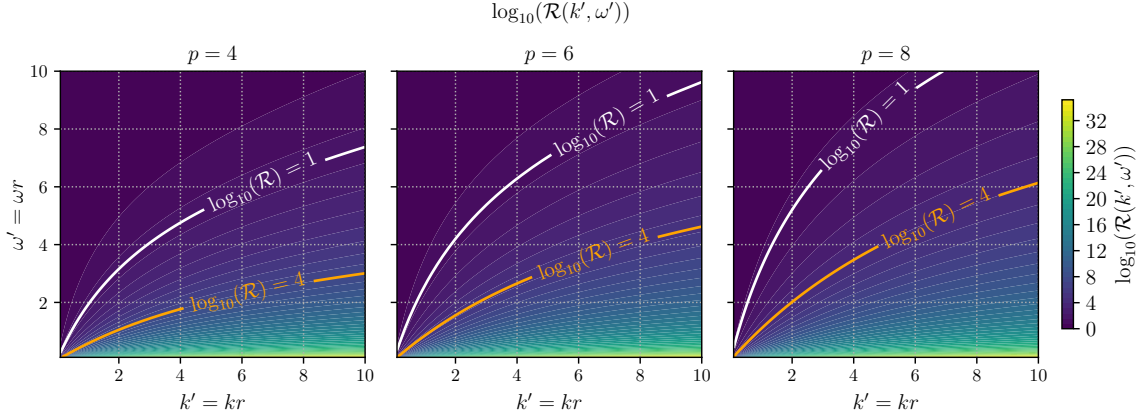}%
    } \caption{ Plot of $\log_{10}(\mathcal{R}(k',\omega'))$ for expansion
    orders $p=4,6,8$. Contour lines indicate $\log_{10}(\mathcal{R}) = 1$
    (white) and $\log_{10}(\mathcal{R}) = 4$ (orange).} \label{fig:gamma_ratio}
\end{figure}

Denoting $k':= kr$ and $\omega':= \omega r$, the truncation error ratio between
\eqref{eq:qbx-flat-panel-exact-truncation-error} and \eqref{eq:qbmax-flat-panel-truncation-error} can be represented as
\begin{equation}
	\mathcal{R}(k',\omega'):= \frac{\gamma(p+1, \sqrt{k'^{2}+ \omega'^{2}} )}{\gamma(p+1, \sqrt{k'^{2}+ \omega'^{2}} - k')}. \label{eq:gamma-ratio}
\end{equation}
\Cref{fig:gamma_ratio} shows the error ratio for $p \in \{ 4,6,8 \}$, along with
level sets of $\log_{10}(\mathcal{R}(k',\omega'))=1$ and
$\log_{10}(\mathcal{R}(k',\omega'))=4$. For near-field evaluation problems, the
expansion radius $r$ at each target point is typically small ($r < 0.1$), which
motivates our choice of plotting domain $[0,10] \times [0,10]$.
The level set $\log_{10}(\mathcal{R})=1$ marks a tenfold reduction
in truncation error from QBX to QBMAX.
For fixed $p$ and $\omega'$, as $k'\to\infty$, the relative QBMAX
truncation error tends to $0$, while the relative QBX truncation error
remains $\mathcal{O}(1)$. While the ratio is derived for the infinite-flat-boundary
model, we demonstrate numerically in \Cref{sec:numerical-results} that this
relationship continues to hold well in more general cases, even for large values of $k$,
particularly for geometries with low-to-moderate curvature relative to
the scale of the expansion.
\subsection{QBMAX Expansion}

We define the new \emph{target-specific} line expansion -- \textit{Quadrature by
Matched Asymptotic Expansion} (QBMAX) as follows.
\begin{definition}[QBMAX expansion]
Suppose $u$ is a layer potential and $\b x \in \partial \Omega$ with an expansion
line $\b l(\tau)$. We define the target-specific asymptotic scaling function
$u_{\text{ma}}$ as $u_{\text{ma}}(\tau) = \exp(k|\b l(\tau)-\b x|)$ for $\tau \in
[0,1]$. Denoting the $p^{\text{th}}$-order matched asymptotic expansion operator
by $\mathcal{T}_p^{\textup{ma}}$, the QBMAX line expansion of $u \circ \b l$ centered at
$\tau = 0$ is defined as
\begin{equation}
\begin{split}
\mathcal{T}_p^{\textup{ma}}[u \circ \b l](\tau) &:= \frac{\mathcal{T}_{p}[(u\circ \b l) \cdot u_{\text{ma}}](\tau)}{u_{\text{ma}}(\tau)} \\
&= \exp(-k r (1-\tau)) \sum_{m=0}^{p} \frac{\partial^{m}_\tau [(u \circ \b l)(\tau) \cdot \exp(k r(1-\tau))] }{m!}\bigg\rvert_{\tau=0} \tau^{m} \\
&= \exp(k r\tau) \sum_{m=0}^{p} \frac{\partial^{m}_\tau [(u \circ \b l)(\tau) \cdot \exp(-k r\tau)] }{m!}\bigg\rvert_{\tau=0} \tau^{m}.
\end{split}
\label{eq:QBMAX}
\end{equation}
\end{definition}
Compared to the previous QBX line expansion \eqref{eq:line-qbx-taylor-general}, the
QBMAX method includes the modification term $u_{\text{ma}}$, which compensates
for the exponential decay of the layer potential along the expansion line.
To recover the original layer potentials,
$u_{\text{ma}}(\tau)$ is needed in the denominator of the new
expansion since the Taylor expansion is computed for $(u \circ \b l) \cdot
u_{\text{ma}}$ rather than $u \circ \b l$ directly.

Direct application of Taylor's theorem (with integral remainder) shows that
\begin{equation}
\begin{split}
    (u \circ \b l) (\tau)  - \mathcal{T}_p^{\textup{ma}}[u \circ \b l](\tau) &=
    e^{kr\tau}\left[e^{-kr\tau} \cdot  (u \circ \b l)(\tau) - \mathcal{T}_p[(u\circ \b l)(\tau) \cdot e^{-kr\tau}](\tau)\right] \\
    &= e^{kr\tau} \left[\frac{\tau^{p+1}}{p!} \int_0^1 (1-\xi)^p \partial_s^{p+1} [(u \circ \b l)(s) \cdot e^{-krs}]\big\rvert_{s = \xi \tau} d\xi\right]
\end{split}
\label{eq:QBMAX-truncation-error}
\end{equation}

\begin{lemma}
   For $m \in \N$,
\begin{align*}
    \partial_\tau^m [(u \circ \b l)(\tau) \cdot e^{-kr\tau}] = e^{-kr\tau} (\partial_\tau - kr)^m (u \circ \b l)(\tau).
\end{align*}
\end{lemma}
\begin{myproof}
Denote $f(\tau) = (u \circ \b l)(\tau)$, $g(\tau) = e^{-kr\tau}$ and $M_g$ be the multiplication operator
defined as $M_g f = g f$. Then,
\begin{align*}
    \partial_{\tau} (fg) = (\partial_\tau f) g + f (\partial_\tau g) = (\partial_\tau f) g - kr f g = g (\partial_\tau - kr) f = M_g (\partial_\tau - kr) f.
\end{align*}
Since $g$ is an eigenfunction of $\partial_{\tau}$ with eigenvalue $-kr$,
$\partial_{\tau} \circ M_g = M_g \partial_{\tau} - kr M_g = M_g (\partial_\tau - kr)$. Hence,
\begin{align*}
\partial_{\tau}^{m} (fg) = \partial_{\tau}^{m-1} \circ M_g (\partial_\tau - kr) f = M_g (\partial_\tau - kr)^m f.
\end{align*}
\end{myproof}
\noindent
Then, the truncation error of the QBMAX expansion in \eqref{eq:QBMAX-truncation-error} becomes
\begin{align*}
    \left|(u \circ \b l)(\tau)  - \mathcal{T}_p^{\textup{ma}}[u \circ \b l](\tau)\right| &=
    e^{kr\tau} \left|\frac{\tau^{p+1}}{p!} \int_0^1 (1-\xi)^p e^{-kr\xi\tau} (\partial_s - kr)^{p+1} (u \circ \b l)(s)\big\rvert_{s = \xi \tau} d\xi\right| \\
    &\leq \frac{\tau^{p+1}}{p!} \int_0^1 (1-\xi)^p e^{kr\tau(1-\xi)} \left|(\partial_s - kr)^{p+1} (u \circ \b l)(s)\big\rvert_{s = \xi \tau}\right| d\xi\\
    &\leq \frac{(r\tau)^{p+1}}{(p+1)!} \sup_{s \in [0,\tau]} e^{kr(\tau - s)}\left|\big(\b n(\b x)\cdot \nabla_{\b z} - k\big)^{p+1} u(\b z)\big\rvert_{\b z = \b l(s)}\right|.
\end{align*}
\begin{remark}
The asymptotic convergence rate of the QBMAX expansion is the same as that of the QBX expansion \eqref{eq:QBX-truncation}, namely $(r\tau)^{p+1}$.
The difference is in the truncation constant.
For modified Helmholtz layer potentials with large $k$, the leading-order decay in the expansion direction is $\exp(-kr(1-\tau))$ 
(this can be shown by a local-coordinate Taylor expansion of the geometry, density, and kernel). 
QBMAX removes this factor before Taylor approximation. The infinite-flat-boundary analysis makes the reduction explicit, 
and the numerical results in \Cref{sec:numerical-results} show the improvement for large $k$ in practice.
\end{remark}

\subsection{QBMAX Quadrature Error}
\label{sec:quad-error}
Suppose $\b x \in \partial \Omega$ and $\b l(\tau)$ is its expansion
line so that $\b l(0) = \b c_{\b x}$ and $\b l(1) = \b x$. The $p^{\text{th}}$-order QBX and QBMAX expansion kernels are defined by
\begin{align}
		\mathcal{K}_{\textup{QBX}}^{p}(\b l(\tau),\b y) &:=
		    \sum_{m=0}^{p} \frac{\partial_{\tau}^{m}
		        \left[ \tau \mapsto \mathcal{K}(\b l(\tau),\b y) \right]}{m!}
		    \bigg\rvert_{\tau=0} \tau^{m}, \label{eq:qbx-kern} \\
	\begin{split}
		\mathcal{K}_{\textup{QBMAX}}^{p}(\b l(\tau),\b y) &:=
		    \exp(-k|\b l(\tau)-\b x|)
		    \sum_{m=0}^{p} \frac{\partial_{\tau}^{m}
		        \left[ \tau \mapsto \mathcal{K}(\b l(\tau),\b y) \cdot \exp(k|\b l(\tau)-\b x|) \right]}{m!}
		    \bigg\rvert_{\tau=0} \tau^{m} \\
	        &= \exp(kr\tau) \sum_{m=0}^{p} \frac{\partial_{\tau}^{m}
	            \left[ \tau \mapsto \mathcal{K}(\b l(\tau),\b y) \cdot \exp(-kr\tau) \right]}{m!}
	            \bigg\rvert_{\tau=0} \tau^{m},
	\end{split} \label{eq:qbmax-kern}
\end{align}
where $\mathcal{K}$ is the kernel function and $r:=|\b l(0) - \b x|$ is the expansion radius at $\b x$. By Leibniz's rule,
\begin{align*}
    \partial_{\tau}^{m}\!\left[\mathcal{K}(\b l(\tau),\b y)\cdot\exp(-kr\tau)\right]\bigg\rvert_{\tau=0}
    = \sum_{j=0}^{m} \binom{m}{j} (-kr)^{m-j}\,
        \partial_{\tau}^{j}\mathcal{K}(\b l(\tau),\b y)\big\rvert_{\tau=0}.
\end{align*}
As a consequence, $\mathcal{K}_{\textup{QBMAX}}^{p}(\b l(\tau),\b y)$ can be rewritten as
\begin{align*}
    \mathcal{K}_{\textup{QBMAX}}^{p}(\b l(\tau),\b y)
    &= \exp(kr\tau) \sum_{m=0}^{p} \sum_{j=0}^{m} \binom{m}{j} \frac{(-kr)^{m-j}}{m!} \partial_{\tau}^{j}\mathcal{K}(\b l(\tau),\b y)\big\rvert_{\tau=0} \tau^{m}\\
    &= \exp(kr\tau) \sum_{j=0}^{p} \sum_{m=j}^{p} \binom{m}{j} \frac{(-kr)^{m-j}}{m!} \partial_{\tau}^{j}\mathcal{K}(\b l(\tau),\b y)\big\rvert_{\tau=0} \tau^{m}\\
    &= \sum_{j=0}^{p} B_{j,p}(kr,\tau)\,
        \partial_{\tau}^{j}\mathcal{K}(\b l(\tau),\b y)\big\rvert_{\tau=0},
\end{align*}
where
\begin{align*}
    B_{j,p}(kr,\tau) := \exp(kr\tau) \sum_{m=j}^{p} \binom{m}{j}
        \frac{(-kr)^{m-j}\,\tau^{m}}{m!}.
\end{align*}
In \Cref{sec:appendix-quad-error},
we show that the quadrature error of the QBX line expansion is of order $\left(\frac{h}{4r}\right)^{2q}$ for on-surface evaluation of the single-layer potential,
consistent with the original QBX series expansion estimate in \Cref{thm:qbx-series-quadrature-error}.
For fixed $k$, $p$, $r$, and $\tau$, the coefficients $B_{j,p}(kr,\tau)$ are finite. Since $\mathcal{K}_{\textup{QBMAX}}^{p}$ is a finite linear combination of
$\partial_{\tau}^{j}\mathcal{K}(\b l(\tau),\b y)\big\rvert_{\tau=0}$, it follows that
\[
 |\text{Quadrature Error}_{\mathrm{QBMAX}}|
 \leq C(k,p,r,\tau,q,\partial\Omega)
 \left(\frac{h}{4r}\right)^{2q}
 \|\sigma\|_{C^{2q}(\partial\Omega)}.
\]
Thus QBMAX retains the same order in $h$ with these parameters fixed. This
argument does not assert that the constant is uniform as $k$ or $kr$ grows.

\section{Numerical Results}\label{sec:numerical-results}
In this section, we present numerical results to demonstrate the performance of
QBMAX compared to QBX methods. First, we compare the cost of forming line
expansions by QBX and QBMAX methods. Then, we consider the accuracy of the QBMAX
method in both 2D and 3D. The numerical results are obtained using the Python
package Pytential \cite{klocknerInducerPytential2025} and its dependencies. The code is available at
\url{https://github.com/ShawnL00/qbmax-paper-experiments}.

\subsection{Cost Analysis}\label{sec:Cost-Analysis}
\begin{figure}[H]
	\centering \resizebox{1\textwidth}{!}{%
		\input{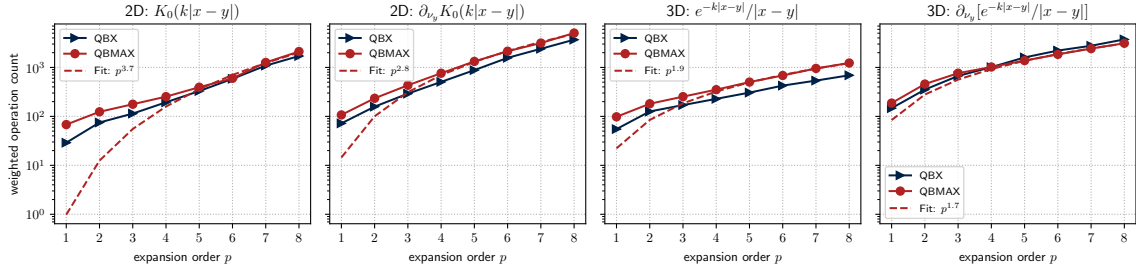}
	} \caption{Weighted symbolic operation-count comparison for 2D/3D expansion
	kernels generated with SymPy. The empirical power-law exponent fitted over
	$5\leq p\leq8$ for QBMAX is reported in the legend.}
	\label{fig:cost}
\end{figure}

We begin by examining the flop count for forming the expansion kernels
$\mathcal{K}_{\textup{QBX}}^{p}$ \eqref{eq:qbx-kern} and
$\mathcal{K}_{\textup{QBMAX}}^{p}$ \eqref{eq:qbmax-kern}. The additional cost of
\eqref{eq:qbmax-kern} relative to \eqref{eq:qbx-kern} comes primarily from the
exponential term and its derivatives. We use the following assumptions for the
flop count:

\vspace{0.5em} \noindent\textbf{Assumptions.}
\begin{enumerate}
    \item The cost of evaluating the exponential function is estimated at $10$
          flops, while other basic function evaluations are counted as $1$ flop.
    \item Derivative computation is performed using the computer algebra package
          SymPy \cite{meurerSymPySymbolicComputing2017a}. After taking each
          of the $p$ derivatives to find the coefficients of the order-$p$ expansion,
          we make use of Sympy's common subexpression elimination (CSE) functionality
          to optimize and reduce the total number of flops.
    \item We provide an estimate of cost scaling with expansion order, via a
          model
          \[
          \text{flop count}\sim p^\alpha.
          \]
          We estimate $\alpha$ by a linear least-squares fit of
          $\log(\text{weighted operation count})$ against $\log(p)$ over
          $p=5,\ldots,8$.
\end{enumerate}
\vspace{0.5em} We report the flop count of both expansions in \Cref{fig:cost},
along with an empirical power-law fit over $5\leq p\leq8$. Under this chosen
weighted symbolic operation-count model, forming the QBMAX kernel expansion has
a cost comparable to that of QBX. The fitted exponent is lower in 3D than in 2D,
which indicates that SymPy's common subexpression elimination is more effective
for the tested 3D expressions over this short range; it should not be interpreted
as an asymptotic complexity result.

Before presenting the numerical results, we introduce a generic problem setup
for evaluating layer potentials.

\subsection{Problem Setup for Layer Potential Evaluation}
Consider a point $\b x \in \partial \Omega$ with an expansion line $\b l(\tau)$. The
goal is to accurately evaluate the regularized layer potential
\begin{equation}
	\int_{\partial \Omega} \mathcal{K}^{p}_{\mathcal{M}}(\b l(\tau),\b y)\sigma(\b y)\,dS_{\b y}, \quad \mathcal{M} \in \{ \text{QBX}, \text{QBMAX} \}, \label{eq:regularized-layer-potential}
\end{equation}
using composite quadrature rules. Specifically, we discretize the boundary
$\partial \Omega$ into $N$ non-overlapping panels $\Gamma_i$ such that
\begin{align*}
\partial \Omega = \bigcup_{i=1}^{N} \Gamma_i.
\end{align*}
On each panel $\Gamma_i$, $q$-point Gauss--Legendre quadrature rules are used
for 2D surfaces and $q^2$-point Gauss--Legendre
tensor-product quadrature rules for 3D surfaces. Although we use tensor-product
quadrature rules for 3D surfaces, we note that the method is equally applicable
to quadrature rules on triangles.

To accurately resolve the kernel, we upsample the quadrature with parameter
$\kappa$. We refer to the nodes before and after upsampling as the
\emph{discretization nodes} and \emph{source nodes}, respectively.

Let $n_{\mathrm{disc}} := q^{n-1}$ and
$n_{\mathrm{up}} := [\kappa(q-1)+1]^{n-1}$ denote the numbers of discretization
and source nodes per panel, respectively, for $n \in \{2,3\}$.
The Nyström method is used to approximate the regularized layer potential
\eqref{eq:regularized-layer-potential}. Let
$\pmb{P} \in \mathbb{R}^{N n_{\mathrm{up}} \times N n_{\mathrm{disc}}}$ be a
block diagonal panelwise interpolation matrix and
$\pmb{W} \in \mathbb{R}^{N n_{\mathrm{up}}\times N n_{\mathrm{up}}}$ be a
diagonal matrix with quadrature weights for the source nodes
$\{\b y_i\}_{i=1}^{N n_{\mathrm{up}}}$. Denote
$\pmb{\sigma}:= [\sigma(\b y_i)]_{i=1, \ldots, N n_{\mathrm{disc}}}$ for the
discretization nodes $\{\b y_i\}$. Then, the discretization of
\eqref{eq:regularized-layer-potential} has the form
\begin{align*}
\sum_{j=1}^{N n_{\mathrm{up}}}\mathcal{K}_{\mathcal{M}}^{p}(\b l(\tau),\b y_j)\pmb{W}_{j,j} [\pmb{P}\pmb{\sigma}]_j.
\end{align*}

We are mainly concerned with $h$-convergence results, where $h$ represents the
mesh size of the panels. In 2D, $h$ is defined as the maximum arc length of the
panels, while in 3D, it is defined as the maximum square root of the area of the
panels. The parameters used in the numerical tests are summarized in
\Cref{tab:parameters}.
\begin{table}[h]
  \centering \caption{Parameters used in the numerical tests.}
  \label{tab:parameters}
  \begin{tabular}{@{}ll@{}}
    \toprule
    \textbf{Parameter} & \textbf{Description} \\
    \midrule
    $p$   & QBX/QBMAX line expansion order \\[2pt]
    $q$   & Quadrature points per panel: $q$ (2D) or $q^{2}$ (3D tensor-product) \\[2pt]
    $\kappa$ & Quadrature upsampling parameter \\[2pt]
    $h$   & Maximum panel size: $\max_{i}\,\mathrm{length}(\Gamma_i)$ (2D) or $\max_{i}\sqrt{\mathrm{area}(\Gamma_i)}$ (3D) \\
    \bottomrule
  \end{tabular}
\end{table}

The experiment-specific discretization parameters are as follows. For the disk
single-layer close-evaluation and convergence studies we use
$(p,q,\kappa)=(5,6,4)$, with $N=40$ for close evaluation and
$N\in\{20,40,80,160,320,480\}$ for convergence. For the corresponding
double-layer studies we use $(p,q,\kappa)=(4,7,4)$, with $N=80$ and
$N\in\{40,80,160,320,640,800\}$, respectively. The disk jump study uses
$(p,q,\kappa)=(6,7,4)$ and
$N\in\{50,100,200,400,600,800\}$. The disk error-ratio study uses
$(p,q,\kappa)=(5,7,4)$, with $N=80$ for the single layer and $N=100$ for
the double layer. In all disk tests, $r=\pi/N$, and the analytic formulas in
\Cref{lem:disk,lem:dis_D} provide the reference values.

For the starfish single- and double-layer close-evaluation studies we use
$(p,q,\kappa)=(6,7,5)$ with $N=80$ and $N=120$, respectively. Their
references use QBMAX with
$(p_{\rm ref},q_{\rm ref},\kappa_{\rm ref})=(9,11,4)$ for the single
layer and $(8,11,4)$ for the double layer
on successively refined meshes up to $N_{\rm ref}=2500$, stopping when
successive values differ in the absolute $\ell^\infty$ norm by less than
$10^{-14}$ for the single layer and $10^{-12}$ for the double layer. The
starfish jump study uses $(p,q,\kappa)=(6,7,4)$ and
$N\in\{100,200,400,600,800,1000\}$. The Helmholtz study uses
$(p,q,\kappa,N)=(6,7,5,300)$. The interior Neumann study uses
$(p,q,\kappa)=(6,7,5)$ and
$N\in\{100,\allowbreak 200,\allowbreak 300,\allowbreak 600,\allowbreak 800,\allowbreak 1000,\allowbreak
1200,\allowbreak 1400,\allowbreak 1600,\allowbreak 1800,\allowbreak 2000\}$; its
reference potential uses $(p_{\rm ref},q_{\rm ref},\kappa_{\rm ref})=(8,11,4)$
and refinement up to $N_{\rm ref}=2500$ with a $10^{-13}$ successive-value
tolerance. The convergence exponents reported for this problem are fitted to
the first five mesh levels.

For the cruller single-layer close-evaluation study we use
$(p,q,\kappa)=(5,6,5)$ and $N=5000$ panels; its reference uses
$(p_{\rm ref},q_{\rm ref},\kappa_{\rm ref})=(6,8,5)$ and $N_{\rm ref}=12800$.
The cruller double-layer jump study uses $(p,q,\kappa)=(4,5,6)$ and
$N\in\{3200,7200,12800,20000,28800\}$; its reported global rates are fitted
to the final four levels. In the starfish and cruller tests, expansion radii
are selected locally by the QBX geometry-refinement procedure; ranges are
reported with the corresponding results below.

We consider the unit disk (\cref{sec:disk}) and starfish domain
(\cref{sec:starfish}) in 2D, and the cruller surface
\cite{barnettHighorderDiscretizationStable2020} (\cref{sec:cruller}) in 3D.
Depending on the test geometry, the reference solution is obtained either
analytically or via self-convergence as noted below. The following test cases
are included for each geometry unless otherwise stated:
\begin{enumerate}
\item \textbf{Close evaluation of layer potentials:} We examine the near-field
      evaluation of layer potentials with analytic densities. For a target point
      $\b x^{*} \in \partial \Omega$ with expansion radius $r_{\b x^{*}}$, we evaluate
      $u(\b x^{*} \mp (1-\tau) r_{\b x^{*}} \b \nu_{\b x^{*}})$ using both QBX and QBMAX
      methods for $\tau \in \{0, 1/8, 1/4, 1/2, 3/4, 7/8, 1\}$ (see
      \Cref{fig:expansion-centers}). Here, $(1-\tau)r_{\b x^{*}}$ represents the distance
      from the evaluation point to the boundary, with $\tau = 1$ corresponding
      to on-surface evaluation and $\tau = 0$ to expansion center evaluation.
      At $\tau=1$, the values corresponding to the minus and plus signs are
      understood as the interior and exterior limits, respectively.

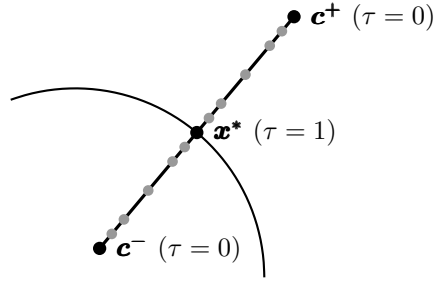
\begin{figure}[H]
  \centering
\begin{tikzpicture}[scale=2.5]
  \draw[domain=0:110,smooth,thick,variable=\x] plot({cos(\x)},{sin(\x)});
  \coordinate (A) at ({cos(50)},{sin(50)});
  \coordinate (B) at ({0.2*cos(50)},{0.2*sin(50)});
  \coordinate (C) at ({1.8*cos(50)},{1.8*sin(50)});
  \draw[very thick] (A) -- (B);
  \draw[very thick] (A) -- (C);
  \foreach \tau in {0,0.125,0.25,0.5,0.75,0.875,1}{
    \fill[gray!80] ($(A)!{1-\tau}!(B)$) circle (0.8pt);
    \fill[gray!80] ($(A)!{1-\tau}!(C)$) circle (0.8pt);
  }
  \fill (A) circle (1pt) node[right=3pt] {$\pmb{\b x^{*}}$ (\(\tau=1\))};
  \fill (B) circle (1pt) node[right=3pt] {$\pmb{\b c^{-}}$ (\(\tau=0\))};
  \fill (C) circle (1pt) node[right=3pt] {$\pmb{\b c^{+}}$ (\(\tau=0\))};
\end{tikzpicture}
  \caption{Illustration of close evaluation along the expansion lines. The
  on-surface target $\b x^*$ is shown with interior expansion center $\b c^{-} = \b x^{*}
  - r_{\b x^{*}}\b \nu_{\b x^{*}}$ and exterior expansion center $\b c^{+} = \b x^{*} +
  r_{\b x^{*}}\b \nu_{\b x^{*}}$. The evaluation points $u(\b x^{*} \mp (1-\tau)
  r_{\b x^{*}}\b \nu_{\b x^{*}})$ for various $\tau$ values are indicated by dots along
  the expansion lines.} \label{fig:expansion-centers}
\end{figure}

\item \textbf{Jump-relation verification:} The jump relations of
      $\mathcal{S}'_{k}\sigma$ and $\mathcal{D}_{k}\sigma$ described in
      \Cref{sec:layer-potentials} are verified numerically. For a set of
      boundary targets $\{\b x_i\}$, we report the achieved accuracy:
\begin{equation}
\begin{split}
&\left\lVert
  \lim_{\substack{\b x \to \b x_i \\ \b x \in \Omega}}
        \mathcal{S}'_{k}\sigma(\b x)
 -\lim_{\substack{\b x \to \b x_i \\ \b x \in \mathbb{R}^{n}\setminus\overline{\Omega}}}
        \mathcal{S}'_{k}\sigma(\b x)
 -\sigma(\b x_i)
\right\rVert_{\ell^{\infty}},
\end{split}
\label{eq:jump-sprime}
\end{equation}

\begin{equation}
\begin{split}
&\left\lVert
  \lim_{\substack{\b x \to \b x_i \\ \b x \in \Omega}}
        \mathcal{D}_{k}\sigma(\b x)
 -\lim_{\substack{\b x \to \b x_i \\ \b x \in \mathbb{R}^{n}\setminus\overline{\Omega}}}
        \mathcal{D}_{k}\sigma(\b x)
 +\sigma(\b x_i)
\right\rVert_{\ell^{\infty}},
\end{split}
\label{eq:jump-d}
\end{equation}
where the $\ell^{\infty}$ norm is taken over all boundary targets $\{\b x_i\}$.
\end{enumerate}

To quantify the order of $h$-convergence, we compute the \emph{estimated order
of convergence} (EOC) using data from successive mesh refinements. Given a
sequence of mesh sizes $\{h_i\}_{i=0}^n$ and their corresponding error measures
$\{e_i\}_{i=0}^n$, the local
$\text{EOC}_i$ between refinement levels is defined as
\begin{align*}
\text{EOC}_i = \frac{\log_{10}(e_i / e_{i-1})}{\log_{10}(h_i / h_{i-1})}, \quad i = 1, \ldots, n.
\end{align*}
To obtain a numerical \textit{global convergence} rate, we perform a first-order
linear regression on the dataset $\{ \left(\log_{10}(h_i),
\log_{10}(e_i)\right)\}_{i=0}^{n}$, minimizing the least squares error.

\subsection{2D Numerical Results -- Forward Evaluation}
\subsubsection{Disk}\label{sec:disk}
The parametrization of the unit disk is given by $(\cos(\theta), \sin(\theta))$
for $\theta \in [0, 2\pi]$. The first example is the close evaluation of the
single layer potential.
\begin{example}[Close evaluation of $\mathcal{S}_{k}\sigma$] \label{ex:disk_S}
This example examines the close evaluation of $\mathcal{S}_{k}\sigma$ with
$\sigma(\b y) = e^{im\theta}$ where $m=7$ and $\theta :=
\operatorname{atan2}(y_2,y_1)$ for $\b y = (y_1,y_2)\in \R^{2}$. By
\Cref{lem:disk},
\begin{align*}
u_{\text{true}}(\b x) := \mathcal{S}_k \sigma (\b x) =\begin{cases}
   e^{im\theta_{\b x}} I_{m}(k|\b x|) K_{m}(k), & |\b x| \leq 1, \\
   e^{im\theta_{\b x}} I_{m}(k) K_{m}(k|\b x|), & |\b x| > 1,
\end{cases}
\end{align*}
where $\theta_{\b x} = \operatorname{atan2}(x_2,x_1)$.

Let $\b x^{*}$ be a discretization node with expansion radius
$r=r_{\b x^{*}}$. For $\tau \in
\{0,\allowbreak 1/8,\allowbreak 1/4,\allowbreak 1/2,\allowbreak 3/4,\allowbreak
7/8,\allowbreak 1\}$, we compute $u(\b x^{*} \mp (1-\tau) r
\b \nu_{\b x^{*}})$ using both QBX and QBMAX methods. The error measure is defined as
\begin{align}
    \text{err}^{\mp}(\tau) := \|u(\b x^{*} \mp (1-\tau) r \b \nu_{\b x^{*}}) - u_{\text{true}}(\b x^{*} \mp (1-\tau)r \b \nu_{\b x^{*}})\|_{\ell^{\infty}}, \label{eq:circle_line_eval_error}
\end{align}
where the $\ell^{\infty}$ norm is taken over all discretization nodes on the
boundary mesh.
\end{example}

\begin{figure}[H]
     \centering
     \begin{subfigure}[t]{0.55\textwidth}
          \centering \resizebox{1.0\textwidth}{!}{%
                \input{numericals_output/circle/S_line_eval_plot.pgf}
          } \caption{} \label{fig:circle_S_direct_eval_a}
     \end{subfigure}\hfill
     \begin{subfigure}[t]{0.4\textwidth}
          \centering
          \includegraphics[width=\textwidth]{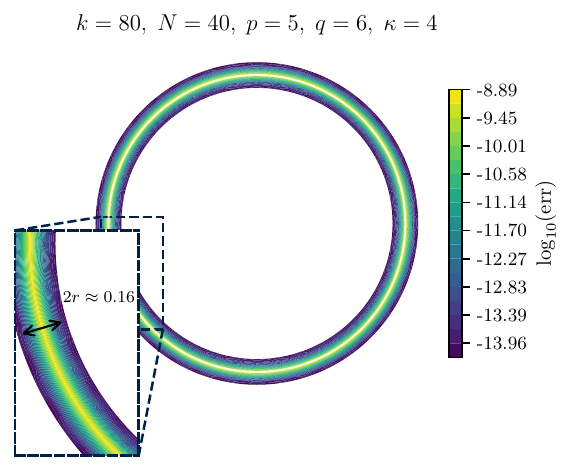}
          \caption{} \label{fig:circle_S_direct_eval_b}
     \end{subfigure}
     \caption{(a) Close evaluation of $u(\b x^{*} \mp (1-\tau) r \b \nu_{\b x^{*}})$ with
     error measure \eqref{eq:circle_line_eval_error}. (b) Visualization of the
     error distribution. The contour width has been exaggerated for enhanced
     visualization.} \label{fig:circle_S_direct_eval}
\end{figure}

In \Cref{fig:circle_S_direct_eval}, the expansion radius $r$ at each
discretization node is approximately $\pi / N \approx 0.0785$. For each value of
$\tau$, the error \eqref{eq:circle_line_eval_error} is computed over $qN = 240$
points. The slope of QBX in \Cref{fig:circle_S_direct_eval_a} is approximately
$6$, consistent with the expected order $p + 1 = 6$. The QBMAX method consistently
outperforms QBX across all values of $k$, maintaining error variations within
two orders of magnitude. \Cref{fig:circle_S_direct_eval_b} visualizes the error
distribution for $k = 80$ (with contour width exaggerated for clarity). At the
expansion centers (outer layer), we achieve $13$-digit accuracy. The accuracy
gradually decreases as we approach the boundary, yielding approximately
$8$-digit accuracy on the boundary itself.

Since the maximum errors in this setting occur at on-surface evaluation ($\tau =
1$), we examine the convergence behavior of both methods under mesh refinements.
\Cref{fig:circle_S_convergence} shows $h$-convergence for on-surface evaluation
of $\mathcal{S}_k\sigma$. The global convergence rate for QBMAX remains robust
at approximately $6$ (interior) and $5.6$ (exterior), with minimal variation
across different $k$ values. In contrast, the QBX method exhibits significant
degradation as $k$ increases, with the global convergence rate declining from
$5.61$ (interior) and $5.46$ (exterior) at $k = 10$ to approximately $3.5$ at $k
= 80$.

\begin{figure}[H]
    \centering
    \begin{subfigure}[t]{0.55\textwidth}
        \vspace{0pt} \centering \resizebox{1.0\textwidth}{!}{%
            \input{numericals_output/circle/S_on_surface_convergence.pgf}
        }
    \end{subfigure}\hfill
    \begin{subfigure}[t]{0.4\textwidth}
        \vspace{0pt} \centering \resizebox{1\textwidth}{!}{%
            \begin{tabular}{|c|c|c|c|c|}
\multicolumn{5}{c}{Global EOC: On-surface Evaluation of $\mathcal{S}_k\sigma$}  \\
\hline
$k$ & {\small QBMAX}$\scriptstyle{(-)}$ & {\small QBX}$\scriptstyle{(-)}$ & {\small QBMAX}$\scriptstyle{(+)}$ & {\small QBX}$\scriptstyle{(+)}$ \\
\hline
10 & 6.26 & 5.61 & 5.64 & 5.46 \\
20 & 6.12 & 5.26 & 5.67 & 5.19 \\
40 & 5.83 & 4.61 & 5.67 & 4.57 \\
80 & 5.87 & 3.56 & 5.60 & 3.53 \\
\hline
\end{tabular}
        }
    \end{subfigure}
    \caption{On-surface evaluation convergence of $\mathcal{S}_{k}\sigma$ for $k
    = 10, 20, 40, 80$.} \label{fig:circle_S_convergence}
\end{figure}

\begin{example}[Close evaluation of $\mathcal{D}_{k}\sigma$] \label{ex:disk_D}
In this example, we consider the close evaluation of $\mathcal{D}_{k}\sigma$
with $\sigma(\b y) = e^{im\theta}$ where $m=7$ and $\theta :=
\operatorname{atan2}(y_2,y_1)$. By \Cref{lem:dis_D},
\[
u_{\text{true}}(\b x) := \mathcal{D}_k \sigma (\b x) =\begin{cases}
  k e^{im \theta_{\b x}} I_{m}(k|\b x|) K'_{m}(k), & |\b x| < 1, \\
  k e^{im \theta_{\b x}} I_{m}'(k) K_{m}(k|\b x|), & |\b x| > 1,
\end{cases}
\]
where $\theta_{\b x} = \operatorname{atan2}(x_2,x_1)$.

Let $\b x^{*}$ be a discretization node with expansion radius
$r=r_{\b x^{*}}$. For $\tau \in
\{0,\allowbreak 1/8,\allowbreak 1/4,\allowbreak 1/2,\allowbreak 3/4,\allowbreak
7/8,\allowbreak 1\}$, we compute $u(\b x^{*} \mp (1-\tau) r
\b \nu_{\b x^{*}})$ using both QBX and QBMAX methods. The error measure is defined as
\begin{align}
   \text{err}^{\mp}(\tau) := \|u(\b x^{*} \mp (1-\tau) r \b \nu_{\b x^{*}}) - u_{\text{true}}(\b x^{*} \mp (1-\tau)r \b \nu_{\b x^{*}})\|_{\ell^{\infty}}, \label{eq:circle_D_error}
\end{align}
where the $\ell^{\infty}$ norm is taken over all discretization nodes on the
boundary mesh.
\end{example}

\begin{figure}[H]
    \centering
    \begin{subfigure}[t]{0.55\textwidth}
        \centering \resizebox{1.0\textwidth}{!}{%
            \input{numericals_output/circle/D_line_eval_plot.pgf}
        } \caption{} \label{fig:circle_D_direct_eval_a}
    \end{subfigure}\hfill
    \begin{subfigure}[t]{0.45\textwidth}
        \centering
        \includegraphics[width=\textwidth]{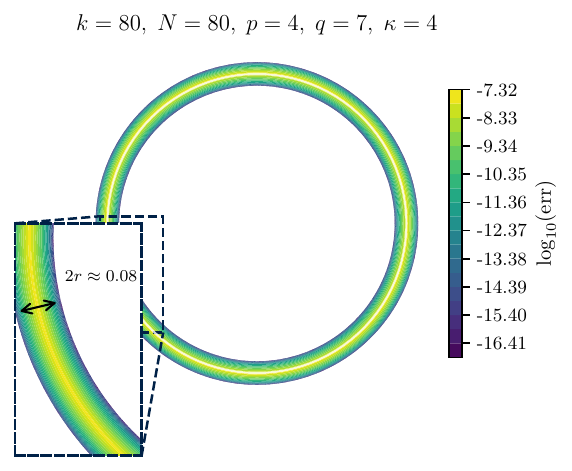}
        \caption{} \label{fig:circle_D_direct_eval_b}
    \end{subfigure}
    \caption{(a) Close evaluation of $u(\b x^{*} \mp (1-\tau) r \b \nu_{\b x^{*}})$ with
    error measure \eqref{eq:circle_D_error}. (b) Visualization of the error
    distribution. The contour width has been exaggerated for enhanced
    visualization.} \label{fig:circle_D_direct_eval}
\end{figure}

\begin{figure}[H]
    \centering
    \begin{subfigure}[t]{0.55\textwidth}
        \vspace{0pt} \centering \resizebox{1.0\textwidth}{!}{%
            \input{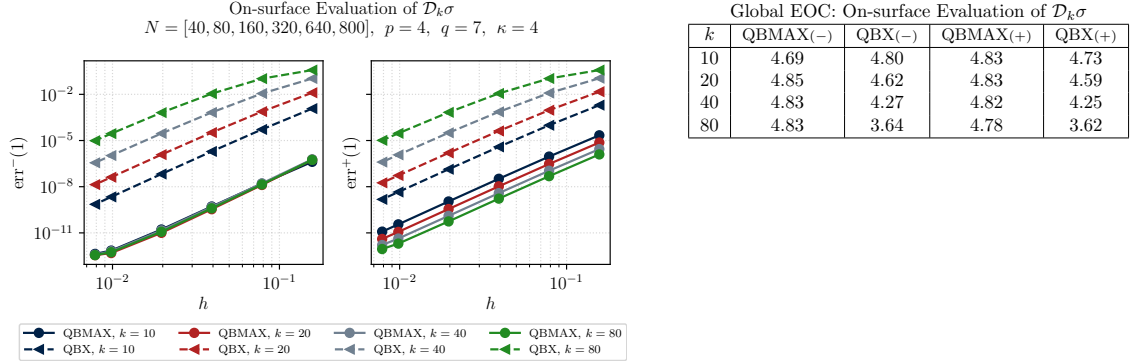}
        }
    \end{subfigure}\hfill
    \begin{subfigure}[t]{0.4\textwidth}
        \vspace{0pt} \centering \resizebox{1\textwidth}{!}{%
            \begin{tabular}{|c|c|c|c|c|}
\multicolumn{5}{c}{Global EOC: On-surface Evaluation of $\mathcal{D}_k\sigma$}  \\
\hline
$k$ & {\small QBMAX}$\scriptstyle{(-)}$ & {\small QBX}$\scriptstyle{(-)}$ & {\small QBMAX}$\scriptstyle{(+)}$ & {\small QBX}$\scriptstyle{(+)}$ \\
\hline
10 & 4.69 & 4.80 & 4.83 & 4.73 \\
20 & 4.85 & 4.62 & 4.83 & 4.59 \\
40 & 4.83 & 4.27 & 4.82 & 4.25 \\
80 & 4.83 & 3.64 & 4.78 & 3.62 \\
\hline
\end{tabular}
        }
    \end{subfigure}
    \caption{On-surface evaluation convergence of $\mathcal{D}_{k}\sigma$ for $k
    = 10, 20, 40, 80$.} \label{fig:circle_D_convergence}
\end{figure}

In \Cref{fig:circle_D_direct_eval}, the expansion radius $r$ at each
discretization node is approximately $\pi / N \approx 0.039$. For each value of
$\tau$, the error \eqref{eq:circle_D_error} is computed over $qN = 560$ points.
We observe that the interior expansion in \Cref{fig:circle_D_direct_eval_a} is
more stable than the exterior expansion, with the convergence rate remaining
close to $5$ for QBMAX. Similarly, \Cref{fig:circle_D_direct_eval_b} provides a
visualization of the error for $k = 80$, and QBMAX achieves $16$-digit accuracy
at expansion centers and $7$-digit accuracy on the boundary.
Figure~\ref{fig:circle_D_convergence} presents the $h$-convergence for
on-surface evaluation of $\mathcal{D}_k\sigma$. The global convergence rate for
QBMAX remains approximately $4.8$ with minimal variation across different $k$
values. On the contrary, the global convergence rate of QBX declines from $4.8$
at $k = 10$ to just $3.6$ at $k = 80$.

\begin{figure}[H]
    \centering
    \begin{subfigure}[t]{0.55\textwidth}
        \vspace{0pt} \centering \resizebox{1.0\textwidth}{!}{%
            \input{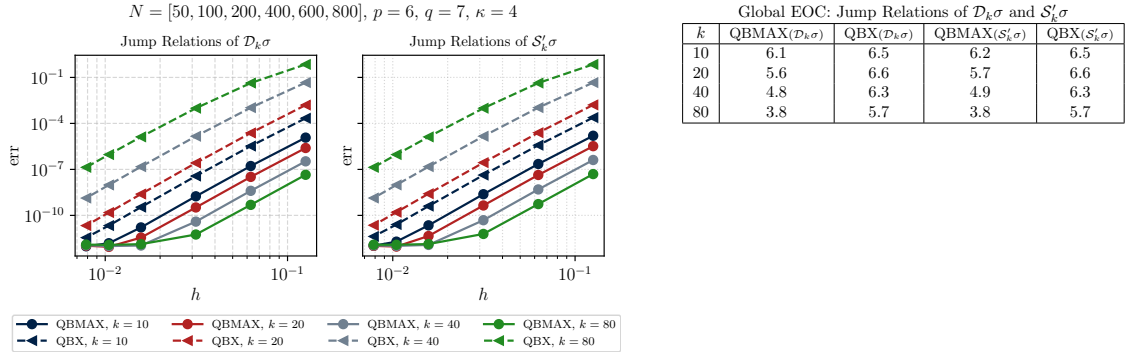}
        }
    \end{subfigure}\hfill
    \begin{subfigure}[t]{0.4\textwidth}
        \vspace{0pt} \centering \resizebox{1\textwidth}{!}{%
        {\normalsize \begin{tabular}{|c|c|c|c|c|}
\multicolumn{5}{c}{\text{Global EOC: Jump Relations of $\mathcal{D}_k\sigma$ and $\mathcal{S}_k'\sigma$}} \\
\hline
$k$ & {\small QBMAX}$\scriptstyle{(\mathcal{D}_k\sigma)}$ & {\small QBX}$\scriptstyle{(\mathcal{D}_k\sigma)}$ & {\small QBMAX}$\scriptstyle{(\mathcal{S}_k'\sigma)}$ & {\small QBX}$\scriptstyle{(\mathcal{S}_k'\sigma)}$ \\
\hline
10 & 6.1 & 6.5 & 6.2 & 6.5 \\
20 & 5.6 & 6.6 & 5.7 & 6.6 \\
40 & 4.8 & 6.3 & 4.9 & 6.3 \\
80 & 3.8 & 5.7 & 3.8 & 5.7 \\
\hline
\end{tabular}}
        }
    \end{subfigure}
   \caption{Convergence of the jump relations of $\mathcal{D}_{k}\sigma$ and
   $\mathcal{S}'_{k}\sigma$ for $k = 10, 20, 40, 80$.}
   \label{fig:circle_jump_relations}
\end{figure}

\begin{example}[Jump relations of $\mathcal{S}'_{k}\sigma$ and $\mathcal{D}_{k}\sigma$]
Define $\sigma(\b y):= \sin(4 \pi y_1) \cos(2\pi y_2)\cos(10 y_2) + 2$ for $\b y=(y_1,
y_2) \in \mathbb{R}^{2}$. We verify the jump relations for
$\mathcal{S}'_{k}\sigma$ and $\mathcal{D}_{k}\sigma$ as in \eqref{eq:jump-sprime}
and \eqref{eq:jump-d}, respectively. The error measures are computed using the
$\ell^{\infty}$ norm over all discretization points at each $h$-refinement
level.
\end{example}
\Cref{fig:circle_jump_relations} shows the convergence of the jump relations for
$\mathcal{D}_{k}\sigma$ and $\mathcal{S}'_{k}\sigma$ for $k = 10, 20, 40, 80$.
We notice that the global EOC is decreasing as $k$ increases for QBMAX. This is because
the QBMAX errors plateau around 12 digits. This error stagnation is likely
caused by floating-point precision limitations when computing small distances
between targets and expansion centers. Nevertheless, examining the convergence
prior to this plateau confirms that the rates remain above 6.

We observe that for small $k$ values, the performance gap between QBX and QBMAX
is not as pronounced as in previous examples. At $k=10$, QBMAX achieves
approximately one additional digit of accuracy. But for each doubling of $k$
after that, QBMAX consistently gains two extra digits in the reported results.

\begin{figure}[H]
    \centering
    \begin{subfigure}[t]{1\textwidth}
        \centering \resizebox{1.0\textwidth}{!}{%
            \input{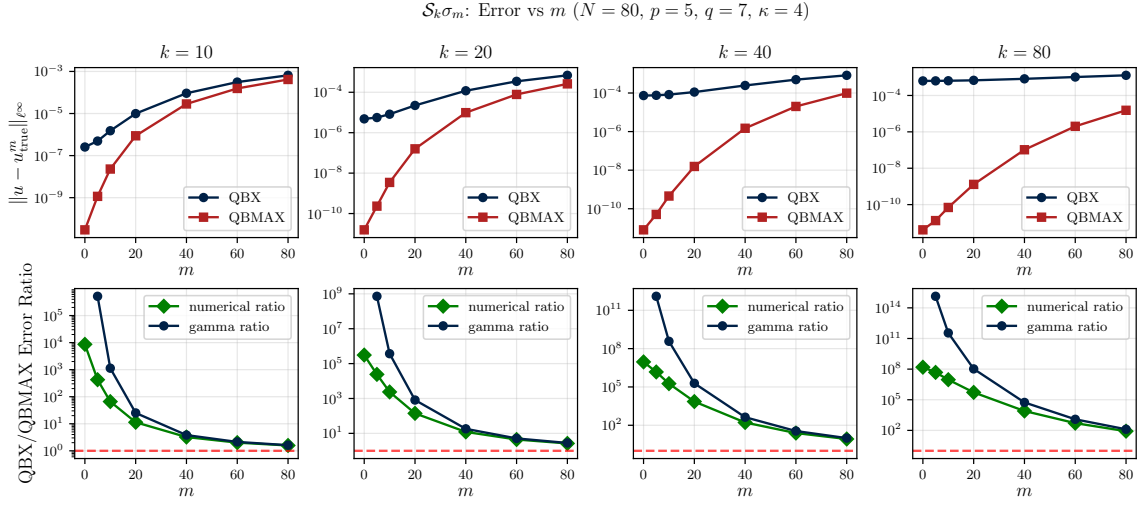}
        } \caption{$\mathcal{S}_k\sigma_m$} \label{fig:circle_S_ratio}
    \end{subfigure}

    \vspace{1em}

    \begin{subfigure}[t]{1\textwidth}
        \centering \resizebox{1.0\textwidth}{!}{%
            \input{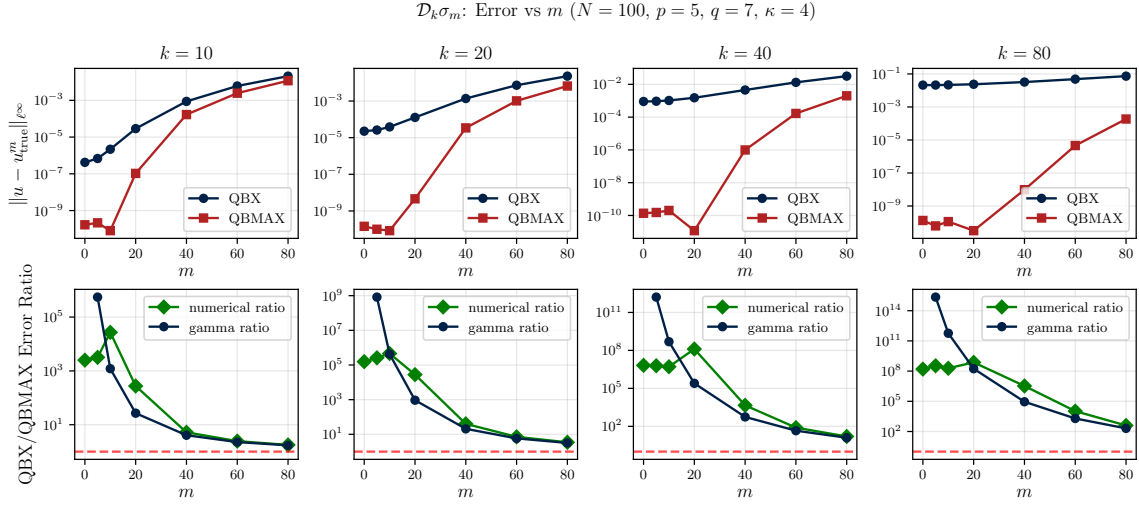}
        } \caption{$\mathcal{D}_k\sigma_m$} \label{fig:circle_D_ratio}
    \end{subfigure}
    \caption{Truncation error comparison for QBX and QBMAX methods with
    oscillatory density $\sigma_{m}$. \textbf{Top row:} $\ell^{\infty}$ error
    over all discretization points on the boundary.
    \textbf{Bottom row:} Ratio of QBX error to QBMAX error based on the data
    from the top row. The gamma ratio from \eqref{eq:gamma-ratio} is included as a
    reference.} \label{fig:circle_ratios}
\end{figure}

Earlier in \Cref{sec:flat-panel-truncation-error}, we derived the truncation
error ratio \eqref{eq:gamma-ratio} for the QBX and QBMAX methods on the flat
panel. Since the disk is a low-curvature geometry, we expect this truncation error
ratio to be representative. Consider the oscillatory density $\sigma_m(\b y) = e^{im
\theta}$ for $m\in \mathbb{Z}$ and $\theta:=\operatorname{atan2}(y_2, y_1)$.
Let $u_{\text{true}}^{m}$ denote the true solution, which is either
$\mathcal{S}_{k}\sigma_m$ or $\mathcal{D}_{k}\sigma_m$. We compute the
$\ell^{\infty}$ error over all discretization points on the boundary and report
results in the top rows of \Cref{fig:circle_ratios}. The performance gap between
QBX and QBMAX decreases as $m$ increases.

The bottom row of \Cref{fig:circle_ratios} shows the ratio of QBX error to QBMAX
error based on the data from the top row. We include the truncation error ratio
estimation from \eqref{eq:gamma-ratio} as a reference. For $m \geq 40$, the gamma
ratio provides a good approximation with variation within roughly one order of
magnitude. In the low-frequency regime ($m \ll k$), the gamma ratio diverges
because QBMAX has zero truncation error on flat boundaries with constant
density.

\subsubsection{Starfish}\label{sec:starfish}
\begin{figure}[htbp]
    \centering
    \resizebox{0.4\textwidth}{!}{\input{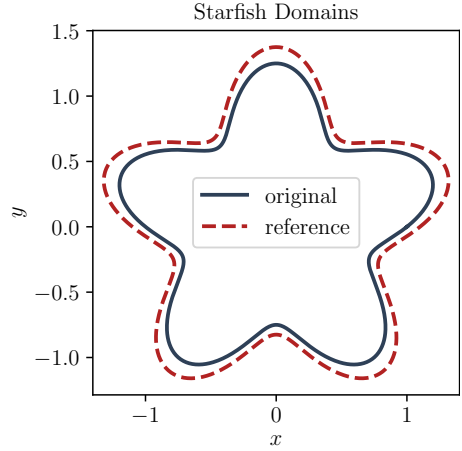}}
    \caption{Starfish domain (solid line) with parametrization $x(\theta) =
    r(\theta)\cos(\theta)$, $y(\theta) = r(\theta)\sin(\theta)$, where $r(\theta) =
    1+\frac{1}{4}\sin(5 \theta)$. The reference domain (dashed line) is scaled:
    $\tilde{x}(\theta) = 1.1\,x(\theta)$, $\tilde{y}(\theta) = 1.1\,y(\theta)$.}
    \label{fig:starfish_domains}
\end{figure}

In this section, we consider the starfish domain as shown in
\Cref{fig:starfish_domains}. The boundary of the starfish domain is
parameterized by $(x(\theta), y(\theta))$, where $x(\theta) = r(\theta)
\cos(\theta)$, $y(\theta) = r(\theta)\sin(\theta)$, and $r(\theta) =
1+\frac{1}{4}\sin(5 \theta)$. For most numerical tests in this section, we
report self-convergence results of the QBX and QBMAX methods since analytic
solutions are generally unavailable. For forward evaluation with analytic
density, the reference solution is obtained by using high-order QBMAX methods on
a refined mesh.

\begin{figure}[H]
    \centering
    \begin{subfigure}[t]{0.55\textwidth}
        \centering \resizebox{1.0\textwidth}{!}{%
            \input{numericals_output/starfish/S_line_eval_plot.pgf}
        } \caption{} \label{fig:starfish_S_direct_eval_a}
    \end{subfigure}\hfill
    \begin{subfigure}[t]{0.4\textwidth}
        \centering
        \includegraphics[width=1.0\textwidth]{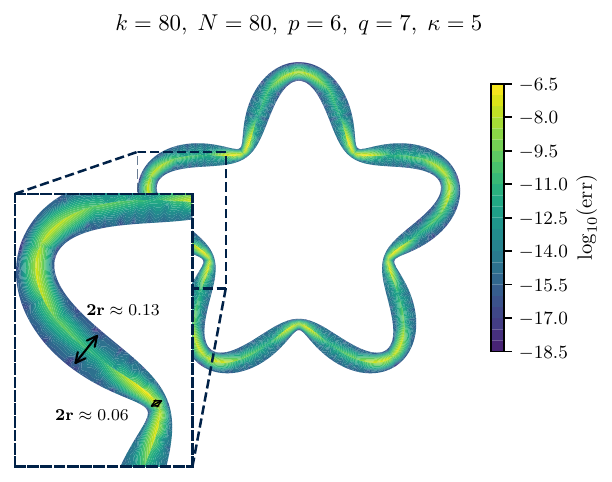}
        \caption{} \label{fig:starfish_S_direct_eval_b}
    \end{subfigure}
   \caption{(a) Close evaluation of $u(\b x^{*} \mp (1-\tau) r_{\b x^{*}} \b \nu_{\b x^{*}})$ with
   error measure \eqref{eq:starfish_line_eval_error}. (b) Visualization of the
   error distribution. The contour width has been exaggerated for enhanced
   visualization.} \label{fig:starfish_S_direct_eval}
\end{figure}

\begin{example}[Close evaluation of $\mathcal{S}_{k}\sigma$] \label{ex:starfish_S}
This example examines the close evaluation of $\mathcal{S}_{k}\sigma$ with
$\sigma(\b y) = \cos(5 \theta) \sin(2 \theta)$ where $\theta :=
\operatorname{atan2}(y_2,y_1)$ for $\b y = (y_1,y_2)\in \R^{2}$.

Let $\b x^{*}$ be a discretization node with expansion radius $r_{\b x^{*}}$. For
$\tau \in \{0, 1/8, 1/4, 1/2, 3/4, 7/8, 1\}$, we compute $u(\b x^{*} \mp (1-\tau)
r_{\b x^{*}} \b \nu_{\b x^{*}})$ using both QBX and QBMAX methods. The error measure is
defined as
\begin{align}
    \text{err}^{\mp}(\tau) := \|u(\b x^{*} \mp (1-\tau) r_{\b x^{*}} \b \nu_{\b x^{*}}) - u_{\text{true}}(\b x^{*} \mp (1-\tau)r_{\b x^{*}} \b \nu_{\b x^{*}})\|_{\ell^{\infty}}, \label{eq:starfish_line_eval_error}
\end{align}
where the $\ell^{\infty}$ norm is taken over all discretization nodes on the
boundary mesh.
\end{example}
For each value of $\tau$, the error \eqref{eq:starfish_line_eval_error} in
\Cref{fig:starfish_S_direct_eval_a} is computed over $qN = 560$ points. The
slopes of both methods are approximately $7$. However, for small $k$ like
$k=10$, there is no significant difference between the two methods, with at
most $1$ order of magnitude of difference.

\begin{example}[Close evaluation of $\mathcal{D}_{k}\sigma$] \label{ex:starfish_D}
In this example, we consider the close evaluation of $\mathcal{D}_{k}\sigma$
with $\sigma(\b y) = \cos(5 \theta) \sin(2 \theta)$ where $\theta :=
\operatorname{atan2}(y_2,y_1)$ for $\b y = (y_1,y_2)\in \R^{2}$.

Let $\b x^{*}$ be a discretization node with expansion radius $r_{\b x^{*}}$. For
$\tau \in \{0, 1/8, 1/4, 1/2, 3/4, 7/8, 1\}$, we compute $u(\b x^{*} \mp
(1-\tau) r_{\b x^{*}} \b \nu_{\b x^{*}})$ using both QBX and QBMAX methods. The error
measure is defined as
\begin{align}
    \text{err}^{\mp}(\tau) := \|u(\b x^{*} \mp (1-\tau) r_{\b x^{*}} \b \nu_{\b x^{*}}) -
    u_{\text{true}}(\b x^{*} \mp (1-\tau)r_{\b x^{*}} \b \nu_{\b x^{*}})\|_{\ell^{\infty}}, \label{eq:starfish_D_line_eval_error}
\end{align}
where the $\ell^{\infty}$ norm is taken over all discretization nodes on the
boundary mesh.
\end{example}
\begin{figure}[H]
    \centering
    \begin{minipage}[t]{0.55\textwidth}
        \vspace{0pt} \centering \resizebox{1.0\textwidth}{!}{%
            \input{numericals_output/starfish/D_line_eval_plot.pgf}
        }
    \end{minipage}\hfill
    \begin{minipage}[t]{0.4\textwidth}
        \vspace{0pt} \centering
        \includegraphics[width=1.0\textwidth]{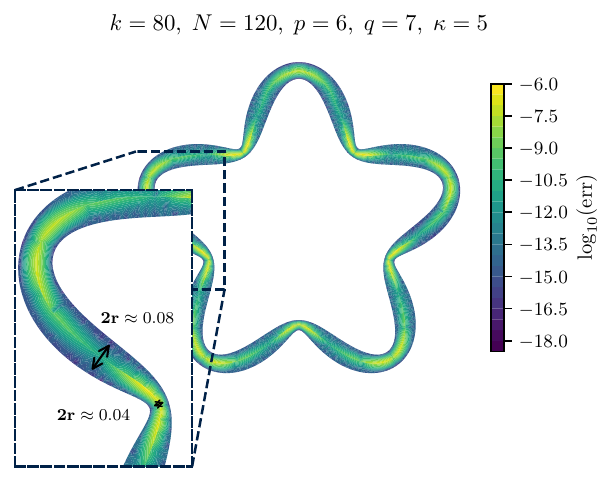}
    \end{minipage}
    \caption{(Left) Close evaluation of $u(\b x^{*} \mp (1-\tau) r_{\b x^{*}} \b \nu_{\b x^{*}})$
    with error measure \eqref{eq:starfish_D_line_eval_error}. (Right)
    Visualization of the error distribution. The contour width has been
    exaggerated for enhanced visualization.} \label{fig:starfish_D_direct_eval}
\end{figure}
Similar to the single layer potential case, \Cref{fig:starfish_D_direct_eval}
shows that the close evaluation of $\mathcal{D}_{k}\sigma$ does not have a
significant difference between the QBX and QBMAX methods for $k=10$.

\begin{figure}[H]
    \centering
    \begin{minipage}{0.7\textwidth}
        \centering \resizebox{0.95\textwidth}{!}{%
            \input{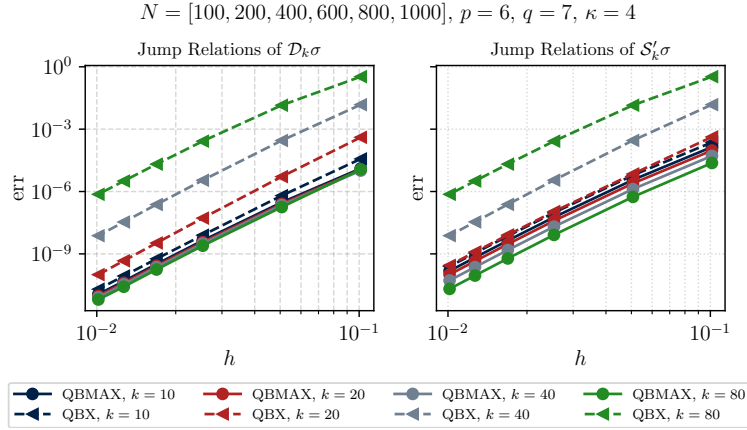}
        }
    \end{minipage}\hfill
    \caption{Convergence of the jump relations of $\mathcal{D}_{k}\sigma$ and
    $\mathcal{S}'_{k}\sigma$ for $k = 10, 20, 40, 80$.}
    \label{fig:starfish_jump_relations}
\end{figure}

\begin{example}[Jump relations of $\mathcal{S}'_{k}$ and $\mathcal{D}_{k}$]
Define $\sigma(\b x):= \sin(4 \pi x_1) \cos(2\pi x_2)\cos(5 x_2) + 2$ for $\b x=(x_1,
x_2)$. We verify the jump relations for $\mathcal{S}'_{k}\sigma$ and
$\mathcal{D}_{k}\sigma$ as in \eqref{eq:jump-sprime} and \eqref{eq:jump-d},
respectively. The error measures are computed using the $\ell^{\infty}$ norm
over all discretization points at each $h$-refinement level.
\end{example}
\Cref{fig:starfish_jump_relations} shows the convergence of the jump relations
for $\mathcal{D}_{k}\sigma$ and $\mathcal{S}_{k}' \sigma$ for $k=10,20,40,80$
with detailed EOC tables in \Cref{tab:starfish_jump_relations}. Both methods
converge at an order of approximately $6$.
\begin{table}
    \centering
    \begin{subtable}[t]{0.95\textwidth}
        \centering \caption{$k = 10$}
        {\footnotesize\setlength{\tabcolsep}{3pt}%
        \begin{tabular}[t]{|c|cc|cc|cc|cc|}
\hline
 $k = 10$ & \multicolumn{2}{c|}{QBMAX$\scriptstyle{(\mathcal{D}_k \sigma)}$} & \multicolumn{2}{c|}{QBX$\scriptstyle{(\mathcal{D}_k \sigma)}$} & \multicolumn{2}{c|}{QBMAX$\scriptstyle{(\mathcal{S}_k'\sigma)}$} & \multicolumn{2}{c|}{QBX$\scriptstyle{(\mathcal{S}_k'\sigma)}$} \\
\hline
$h$ & Error & EOC & Error & EOC & Error & EOC & Error & EOC \\
\hline
1.0e-01 & 1.3e-05 &   & 3.7e-05 &   & 1.4e-04 &   & 2.3e-04 &   \\
5.1e-02 & 3.0e-07 & 5.4 & 6.3e-07 & 5.9 & 3.6e-06 & 5.2 & 6.1e-06 & 5.3 \\
2.5e-02 & 4.4e-09 & 6.1 & 7.9e-09 & 6.3 & 5.8e-08 & 6.0 & 9.6e-08 & 6.0 \\
1.7e-02 & 3.2e-10 & 6.4 & 5.8e-10 & 6.5 & 4.4e-09 & 6.4 & 7.2e-09 & 6.4 \\
1.3e-02 & 4.8e-11 & 6.6 & 8.7e-11 & 6.6 & 6.6e-10 & 6.5 & 1.1e-09 & 6.5 \\
1.0e-02 & 1.1e-11 & 6.5 & 2.0e-11 & 6.6 & 1.5e-10 & 6.6 & 2.5e-10 & 6.6 \\
\hline
 & & 6.1 & & 6.3 & & 6.0 & & 6.0 \\
\hline
\end{tabular}
\begin{tabular}[t]{|>{\centering\arraybackslash}p{2cm}|>{\centering\arraybackslash}p{2cm}|}
\hline
\multicolumn{2}{|c|}{QBX/QBMAX Error Ratio} \\
\hline
$\mathcal{D}_k\sigma$ & $\mathcal{S}_k'\sigma$ \\
\hline
2.90 & 1.69 \\
2.12 & 1.67 \\
1.81 & 1.66 \\
1.80 & 1.65 \\
1.79 & 1.65 \\
1.76 & 1.65 \\
\hline
\end{tabular}}
    \end{subtable}

    \vspace{1em}

    \begin{subtable}[t]{0.95\textwidth}
        \centering \caption{$k = 20$}
        {\footnotesize\setlength{\tabcolsep}{3pt}%
        \begin{tabular}[t]{|c|cc|cc|cc|cc|}
\hline
 $k = 20$ & \multicolumn{2}{c|}{QBMAX$\scriptstyle{(\mathcal{D}_k \sigma)}$} & \multicolumn{2}{c|}{QBX$\scriptstyle{(\mathcal{D}_k \sigma)}$} & \multicolumn{2}{c|}{QBMAX$\scriptstyle{(\mathcal{S}_k'\sigma)}$} & \multicolumn{2}{c|}{QBX$\scriptstyle{(\mathcal{S}_k'\sigma)}$} \\
\hline
$h$ & Error & EOC & Error & EOC & Error & EOC & Error & EOC \\
\hline
1.0e-01 & 1.1e-05 &   & 4.0e-04 &   & 9.4e-05 &   & 4.2e-04 &   \\
5.1e-02 & 2.5e-07 & 5.5 & 5.2e-06 & 6.3 & 2.5e-06 & 5.3 & 7.0e-06 & 5.9 \\
2.5e-02 & 3.7e-09 & 6.1 & 5.2e-08 & 6.6 & 3.9e-08 & 6.0 & 1.1e-07 & 6.0 \\
1.7e-02 & 2.7e-10 & 6.4 & 3.3e-09 & 6.8 & 2.9e-09 & 6.4 & 8.2e-09 & 6.4 \\
1.3e-02 & 4.0e-11 & 6.6 & 4.7e-10 & 6.8 & 4.5e-10 & 6.5 & 1.2e-09 & 6.6 \\
1.0e-02 & 9.5e-12 & 6.5 & 1.0e-10 & 6.9 & 1.0e-10 & 6.6 & 2.8e-10 & 6.6 \\
\hline
 & & 6.1 & & 6.6 & & 6.0 & & 6.2 \\
\hline
\end{tabular}
\begin{tabular}[t]{|>{\centering\arraybackslash}p{2cm}|>{\centering\arraybackslash}p{2cm}|}
\hline
\multicolumn{2}{|c|}{QBX/QBMAX Error Ratio} \\
\hline
$\mathcal{D}_k\sigma$ & $\mathcal{S}_k'\sigma$ \\
\hline
36.29 & 4.49 \\
20.48 & 2.85 \\
14.28 & 2.80 \\
12.41 & 2.78 \\
11.56 & 2.77 \\
10.55 & 2.76 \\
\hline
\end{tabular}}
    \end{subtable}

    \vspace{1em}

    \begin{subtable}[t]{0.95\textwidth}
        \centering \caption{$k = 40$}
        {\footnotesize\setlength{\tabcolsep}{3pt}%
        \begin{tabular}[t]{|c|cc|cc|cc|cc|}
\hline
 $k = 40$ & \multicolumn{2}{c|}{QBMAX$\scriptstyle{(\mathcal{D}_k \sigma)}$} & \multicolumn{2}{c|}{QBX$\scriptstyle{(\mathcal{D}_k \sigma)}$} & \multicolumn{2}{c|}{QBMAX$\scriptstyle{(\mathcal{S}_k'\sigma)}$} & \multicolumn{2}{c|}{QBX$\scriptstyle{(\mathcal{S}_k'\sigma)}$} \\
\hline
$h$ & Error & EOC & Error & EOC & Error & EOC & Error & EOC \\
\hline
1.0e-01 & 1.0e-05 &   & 1.5e-02 &   & 5.2e-05 &   & 1.5e-02 &   \\
5.1e-02 & 2.2e-07 & 5.5 & 2.9e-04 & 5.7 & 1.3e-06 & 5.3 & 2.9e-04 & 5.7 \\
2.5e-02 & 3.2e-09 & 6.1 & 3.5e-06 & 6.4 & 2.0e-08 & 6.0 & 3.5e-06 & 6.4 \\
1.7e-02 & 2.3e-10 & 6.5 & 2.4e-07 & 6.6 & 1.5e-09 & 6.4 & 2.4e-07 & 6.6 \\
1.3e-02 & 3.4e-11 & 6.6 & 3.4e-08 & 6.7 & 2.3e-10 & 6.6 & 3.4e-08 & 6.7 \\
1.0e-02 & 8.1e-12 & 6.5 & 7.5e-09 & 6.8 & 5.3e-11 & 6.6 & 7.6e-09 & 6.8 \\
\hline
 & & 6.1 & & 6.3 & & 6.0 & & 6.3 \\
\hline
\end{tabular}
\begin{tabular}[t]{|>{\centering\arraybackslash}p{2cm}|>{\centering\arraybackslash}p{2cm}|}
\hline
\multicolumn{2}{|c|}{QBX/QBMAX Error Ratio} \\
\hline
$\mathcal{D}_k\sigma$ & $\mathcal{S}_k'\sigma$ \\
\hline
1.5e+03 & 2.9e+02 \\
1.3e+03 & 2.2e+02 \\
1.1e+03 & 1.7e+02 \\
1.0e+03 & 1.6e+02 \\
1.0e+03 & 1.5e+02 \\
9.3e+02 & 1.4e+02 \\
\hline
\end{tabular}}
    \end{subtable}

    \vspace{1em}

    \begin{subtable}[t]{0.95\textwidth}
        \centering \caption{$k = 80$}
        {\footnotesize\setlength{\tabcolsep}{3pt}%
        \begin{tabular}[t]{|c|cc|cc|cc|cc|}
\hline
 $k = 80$ & \multicolumn{2}{c|}{QBMAX$\scriptstyle{(\mathcal{D}_k \sigma)}$} & \multicolumn{2}{c|}{QBX$\scriptstyle{(\mathcal{D}_k \sigma)}$} & \multicolumn{2}{c|}{QBMAX$\scriptstyle{(\mathcal{S}_k'\sigma)}$} & \multicolumn{2}{c|}{QBX$\scriptstyle{(\mathcal{S}_k'\sigma)}$} \\
\hline
$h$ & Error & EOC & Error & EOC & Error & EOC & Error & EOC \\
\hline
1.0e-01 & 1.2e-05 &   & 3.4e-01 &   & 2.4e-05 &   & 3.4e-01 &   \\
5.1e-02 & 1.8e-07 & 6.0 & 1.4e-02 & 4.6 & 5.5e-07 & 5.4 & 1.4e-02 & 4.6 \\
2.5e-02 & 2.5e-09 & 6.2 & 2.6e-04 & 5.7 & 8.2e-09 & 6.1 & 2.6e-04 & 5.7 \\
1.7e-02 & 1.8e-10 & 6.5 & 2.1e-05 & 6.3 & 6.0e-10 & 6.4 & 2.1e-05 & 6.3 \\
1.3e-02 & 2.6e-11 & 6.6 & 3.2e-06 & 6.5 & 9.1e-11 & 6.6 & 3.2e-06 & 6.5 \\
1.0e-02 & 6.4e-12 & 6.4 & 7.3e-07 & 6.6 & 2.1e-11 & 6.6 & 7.4e-07 & 6.6 \\
\hline
 & & 6.3 & & 5.7 & & 6.1 & & 5.7 \\
\hline
\end{tabular}
\begin{tabular}[t]{|>{\centering\arraybackslash}p{2cm}|>{\centering\arraybackslash}p{2cm}|}
\hline
\multicolumn{2}{|c|}{QBX/QBMAX Error Ratio} \\
\hline
$\mathcal{D}_k\sigma$ & $\mathcal{S}_k'\sigma$ \\
\hline
2.9e+04 & 1.4e+04 \\
7.9e+04 & 2.6e+04 \\
1.1e+05 & 3.2e+04 \\
1.2e+05 & 3.4e+04 \\
1.2e+05 & 3.5e+04 \\
1.2e+05 & 3.5e+04 \\
\hline
\end{tabular}}
    \end{subtable}
    \caption{Estimated order of convergence for jump relations at
    $k = 10, 20, 40, 80$. Global convergence rates appear in the final row of
    each EOC column.}
    \label{tab:starfish_jump_relations}
\end{table}

\begin{remark}
    At small $k$ values, QBMAX yields only marginal improvement over QBX on the
    starfish domain, possibly due to the local curvature. We do not yet have a theoretical
    error analysis for curved geometries. However, we expect that adding local
    curvature correction terms to the expansion could improve performance. This
    will be explored in future work.
\end{remark}

\begin{example}[Close evaluation of Helmholtz layer potentials]\label[example]{ex:helmholtz-starfish}
In this example, we show the applicability of the QBMAX method to the
oscillatory Helmholtz layer potential. For $k \in \C$, the Green's function for
the 2D Helmholtz equation is given by
\[
  \mathcal{G}_{k}(\b x,\b y) = \frac{i}{4} H_{0}^{(1)}(k|\b x-\b y|),
\]
where $H_{0}^{(1)}$ is the zeroth order Hankel function of the first kind and
$\b x,\b y \in \R^{2}$. The reference solution is constructed as a sum of point
sources in the following form:
\[
  u_{\text{ref}}(\b x) := \sum_{j=1}^{10} H_{0}^{(1)}(k|\b x-\b y_j|),
\]
where $\b y_j$ are chosen randomly on the boundary of the scaled reference
domain $\tilde\Omega:=1.1\Omega$ (see \Cref{fig:starfish_domains}). By
Green's representation theorem
\cite{coltonInverseAcousticElectromagnetic2013}, for $\b x \in \Omega$,
\begin{equation}
\int_{\pa \Omega} \frac{\pa \mathcal{G}_{k}(\b x,\b y)}{\pa \b \nu_{\b y}} u_{\text{ref}}(\b y) dS_{\b y} - \int_{\pa \Omega} \mathcal{G}_{k}(\b x,\b y) \frac{\pa u_{\text{ref}}(\b y)}{\pa \b \nu_{\b y}} dS_{\b y} + u_{\text{ref}}(\b x) = 0. \label{eq:helmholtz-qbx-eval}
\end{equation}
In the case of Helmholtz layer potentials, we define the target-specific
regularized function as
$u_{\text{ma}}(\tau)=\exp(\operatorname{Im}(k)|\b l(\tau)-\b x^{*}|)$ for target point
$\b x^{*} \in \partial \Omega$. This example is intended to test the rapidly
evanescent regime with $\operatorname{Im}(k)>0$. Because
\eqref{eq:helmholtz-qbx-eval} is the interior Green representation, we evaluate
at $\b x^{*}-(1-\tau)r_{\b x^{*}}\b\nu_{\b x^{*}}$,
with $\tau=1$ understood as the interior limit.
For each discretization node $\b x^{*}$ with expansion radius
$r_{\b x^{*}}$, we evaluate \eqref{eq:helmholtz-qbx-eval} using both
QBX and QBMAX methods with error measure
\begin{align}
    \text{err}(\tau) := \|u(\b x^{*} - (1-\tau) r_{\b x^{*}} \b \nu_{\b x^{*}}) - u_{\text{ref}}(\b x^{*} - (1-\tau)r_{\b x^{*}} \b \nu_{\b x^{*}})\|_{\ell^{\infty}}, \label{eq:starfish_helmholtz_error}
\end{align}
where the $\ell^{\infty}$ norm is taken over all discretization nodes on the
boundary mesh.
\end{example}

\begin{figure}[H]
    \centering
    \resizebox{0.8\textwidth}{!}{%
        \input{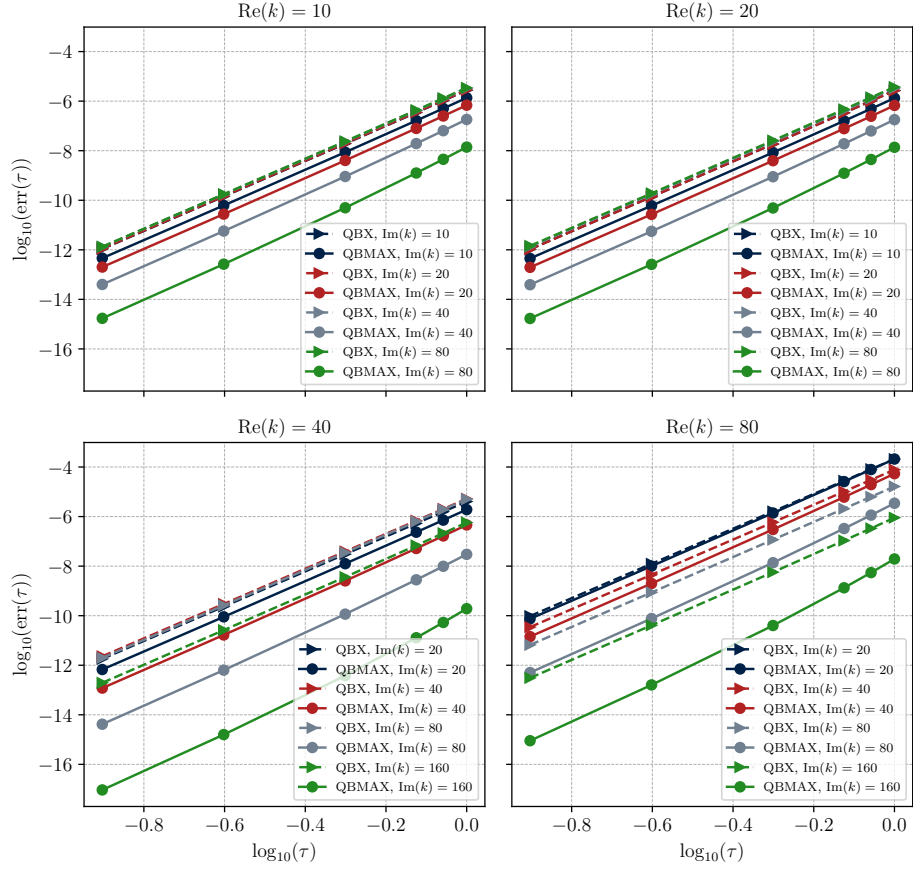}
    }
    \caption{Close evaluation of Helmholtz layer potentials with complex
    wavenumbers over the indicated real and imaginary parts of $k$, using the
    error measure \eqref{eq:starfish_helmholtz_error}.}
    \label{fig:starfish_helmholtz_direct_eval}
\end{figure}

\Cref{fig:starfish_helmholtz_direct_eval} shows the close evaluation of
Helmholtz layer potentials with complex wavenumbers having $\operatorname{Re}(k)
= 10, 20, 40, 80$ and various imaginary components. The error measure
\eqref{eq:starfish_helmholtz_error} is over $qN = 2100$ target points with
expansion radius ranging from $0.008$ to $0.017$. When $\operatorname{Im}(k) \leq
\operatorname{Re}(k)$, there is no significant difference in performance between
the two methods for reported $k$ values. Nevertheless, as $\operatorname{Im}(k) >
\operatorname{Re}(k)$, QBMAX consistently outperforms QBX, often achieving at least
one additional digit of accuracy. This performance gap becomes increasingly
pronounced with larger values of $\operatorname{Im}(k)$.

\subsection{2D Numerical Results -- Interior Neumann Problem}
In this section, we solve an interior Neumann boundary value problem on the
starfish domain as in \Cref{fig:starfish_domains},
\[
\begin{cases}
 -\Delta u + k^{2} u = 0, & \text{in } \Omega, \\
 \partial_{\b \nu} u = f, & \text{on } \partial \Omega.
\end{cases}
\]
To establish a reference solution $u_{\text{ref}}$, we use the scaled reference
domain $\tilde\Omega=1.1\Omega$ as shown in \Cref{fig:starfish_domains}. For a given density
$\tilde{\sigma}$ on $\partial \tilde{\Omega}$, $u_{\text{ref}}$ is defined by
\begin{align}
  u_{\text{ref}}(\b x):= \int_{\partial \tilde{\Omega}} \mathcal{G}_{k}(\b x,\b y) \tilde{\sigma}(\b y) \,dS_{\b y}, \qquad  \b x \in \mathbb{R}^{2}.
  \label{eq:reference_solution}
\end{align}
This reference solution is used to generate the Neumann data on $\partial
\Omega$.

The Neumann problem is represented as a second-kind integral formulation and
discretized with the Nyström method. The resulting linear system is solved
via GMRES (Generalized Minimal Residual Method) to determine the discrete
density on $\partial \Omega$. The GMRES iteration count is recorded with a
relative residual error tolerance set to $10^{-12}$ in the $\ell^{2}$ norm from a
zero initial guess.

To test the accuracy of the solution, suppose $\b \gamma: [0, 1] \to \partial
\Omega$ is the parametrization of the starfish boundary. We generate a set of
1000 target points $\{ \b x_i \}_{i=1}^{1000}$ on $\partial \Omega$ by selecting
equispaced points $t_i = (i-1)/1000$ in $[0, 1)$ and setting $\b x_i =
\b \gamma(t_i)$. For each target point $\b x_i$, we define the local expansion radius
as
\[
r_{i} = \frac{|\b \gamma'(t_i)|}{2N},
\]
where $N$ is the number of panels in the boundary discretization.

The reference solution $u_{\text{ref}}$ is constructed as the single-layer
potential $u_{\text{ref}} = \mathcal{S}_{k}\tilde{\sigma}$ with density
\[
   \tilde{\sigma}(\b y) = \cos(10 \theta) \sin (5 \theta) + 2, \quad \text{where } \theta = \operatorname{atan2}(y_2,y_1).
\]
The corresponding Neumann data $f$ is defined by $f := \partial_{\b \nu}
u_{\text{ref}}|_{\partial \Omega}$. The solution can be represented as $u =
\mathcal{S}_{k}\sigma$ for some unknown density $\sigma$ on $\partial \Omega$.
This leads to the second-kind integral equation
\begin{align}
\frac{1}{2} \sigma (\b x) + \mathcal{S}'_{k}\sigma(\b x) = f(\b x), \quad \b x \in \partial \Omega. \label{eq:neumann-formulation-starfish}
\end{align}
For $\tau \in \{\frac{1}{2}, 1\}$, we report the relative error:
\[
  \text{relative error }(\tau) :=  \frac{ \max_i |u(\b x_i - (1 -\tau) r_i \b \nu_{\b x_i}) - u_{\text{ref}}(\b x_i - (1 - \tau) r_i \b \nu_{\b x_i})|}{\max_i|u_{\text{ref}}(\b x_i - (1 - \tau) r_i \b \nu_{\b x_i})|}.
\]

\begin{figure}[htbp]
   \centering \resizebox{\columnwidth}{!}{%
       \input{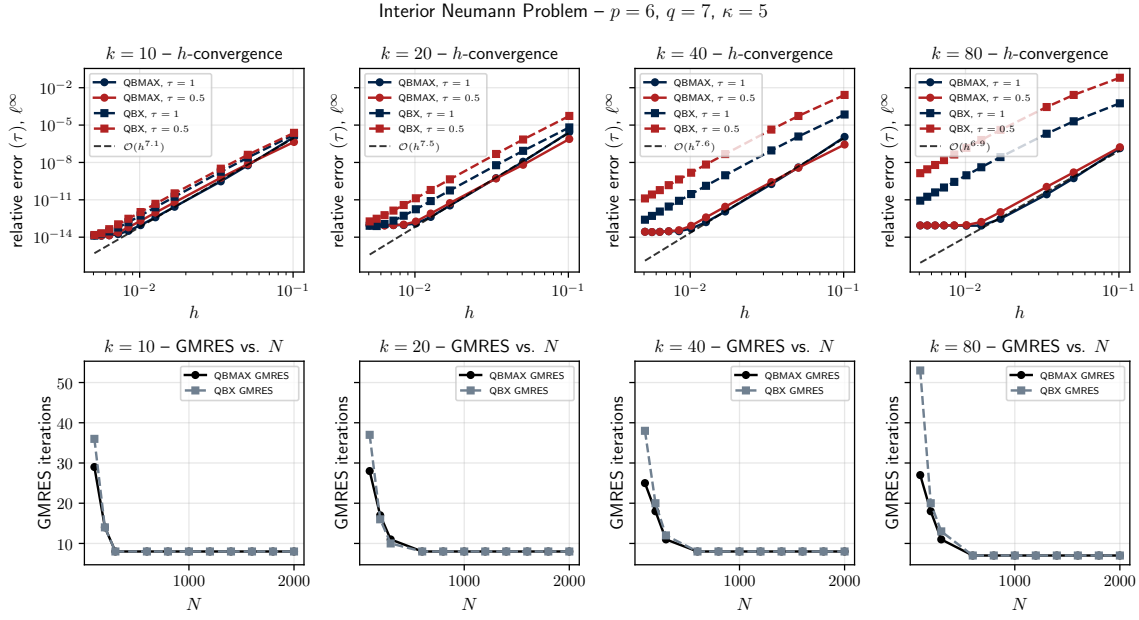}
   } \caption{Convergence results for the interior Neumann problem with
   $h$-refinement for different $k$ values. \textbf{Top row:} Relative error in
   $\ell^\infty$ norm versus element size $h$ for $\tau \in \{1, 0.5\}$.
   \textbf{Bottom row:} GMRES iteration count versus number of boundary
   elements $N$ for both methods with relative tolerance set to $10^{-12}$.}
   \label{fig:starfish_neumann_convergence}
\end{figure}

The top row in \Cref{fig:starfish_neumann_convergence} shows the convergence of
the relative error for different panel sizes $h$ at $\tau = 1$ and $\tau = 0.5$.
The estimated convergence orders shown in the top-row legends range
from about $6.9$ to $7.6$, consistent with the expected order $p+1=7$.
The bottom row shows the number of GMRES iterations required to solve the
linear system with relative residual tolerance set to $10^{-12}$. Outside the
underresolved cases, GMRES reaches the target relative residual in about 7--8
iterations.

\subsection{3D Numerical Results}
\subsubsection{Cruller}\label{sec:cruller}
In this section, we consider the cruller surface
\cite{barnettHighorderDiscretizationStable2020}, which is parameterized by
\begin{equation}
    \b r(u,v) = \begin{pmatrix}
    \cos(u)[R + r H(u,v) \cos(v)] \\
    \sin(u)[R + r H(u,v) \cos(v)] \\
    r H(u,v) \sin(v)
    \end{pmatrix}, \quad (u,v) \in [0,2\pi) \times [0,2\pi), \label{eq:cruller_surface}
\end{equation}
where $R = 4$, $r = 2$, and $H(u,v) = 1 + 0.25 \cos(5u + 3v)$. To construct a
set of target points $\{ \b x_{i,j} \}_{i,j=1}^{100} \subset \partial\Omega$, we
define uniformly spaced parameter values $u_i = 2\pi(i-1) / 100, v_j = 2\pi(j-1)
/ 100 \text{ for } i,j = 1, \ldots, 100, $ and set $\b x_{i,j} := \b r(u_i,
v_j)$. For brevity, we refer to these target points as $\b x^{*}$ in the following
examples.
\begin{example}[Close evaluation of $\mathcal{S}_{k}\sigma$]
We examine the close evaluation of $\mathcal{S}_{k}\sigma$ with density
$\sigma(\b x) = x_1 \sin(2 x_2)\cos(5x_3) + 2$ for $\b x = (x_1,x_2,x_3) \in
\mathbb{R}^{3}$. The reference solution $u_{\text{ref}}$ at each target point
$\b x^{*}$ is generated with the QBMAX method on a refined source geometry.

Let $\b x^{*}$ be a target point with expansion radius $r_{\b x^{*}}$. For $\tau \in
\{0, 1/8, 1/4, 1/2, 3/4, 7/8, 1\}$, we compute $u(\b x^{*} \mp
(1-\tau)r_{\b x^{*}}\b \nu_{\b x^{*}})$ using both QBX and QBMAX methods. The error
measure is defined as
\begin{equation}
    \text{err}^{\mp}(\tau) := \|u(\b x^{*} \mp (1-\tau) r_{\b x^{*}} \b \nu_{\b x^{*}}) - u_{\text{ref}}(\b x^{*} \mp (1-\tau)r_{\b x^{*}}\b \nu_{\b x^{*}})\|_{\ell^{\infty}}, \label{eq:cruller_s_line_eval_error}
\end{equation}
where the $\ell^{\infty}$ norm is taken over all target points $\b x^{*}$ for each
value of $\tau$.
\end{example}

\begin{figure}[H]
    \centering \resizebox{0.7\textwidth}{!}{%
            \input{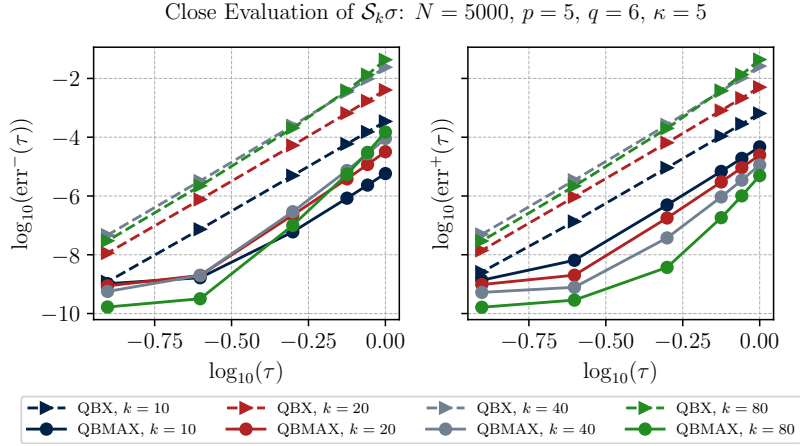}}
            \caption{Close evaluation of $u(\b x^{*} \mp (1-\tau) r_{\b x^{*}}
            \b \nu_{\b x^{*}})$ with error measure \eqref{eq:cruller_s_line_eval_error}.}
            \label{fig:cruller_S_direct_eval}
\end{figure}
Close evaluation results for $\mathcal{S}_{k}\sigma$ along the expansion lines
are presented in \Cref{fig:cruller_S_direct_eval}. The cruller surface is
discretized with $N = 5000$ panels, and the expansion radius ranges between
$0.047$ and $0.11$ depending on the local geometries of the target points. The
slope of both methods is approximately $6$.

\begin{example}[Jump relations of $\mathcal{D}_{k}\sigma$]
Define $\sigma(\b x):= \sin(3 \pi x_1) \cos(5 \pi x_2) \sin(10 x_3) + 2$ for $\b x =
(x_1,x_2,x_3) \in \mathbb{R}^3$. We verify the jump relation \eqref{eq:jump-d} for
$\mathcal{D}_{k}\sigma$, with the $\ell^{\infty}$ norm taken over all
target points $\{\b x^{*}\}$.
\end{example}

\begin{figure}[H]
    \centering
    \begin{subfigure}[t]{0.5\textwidth}
        \centering \resizebox{1.0\textwidth}{!}{%
            \input{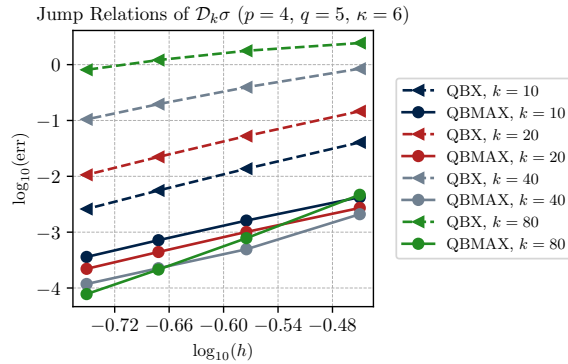}
        }
    \end{subfigure}
    \caption{Jump relation convergence for $\mathcal{D}_{k}\sigma$ on the
    cruller surface for $k = 10, 20, 40, 80$.} \label{fig:cruller_D_convergence}
\end{figure}

\begin{figure}[H]
    \centering
    \begin{subfigure}[t]{0.48\textwidth}
        \centering \includegraphics[width=\textwidth]{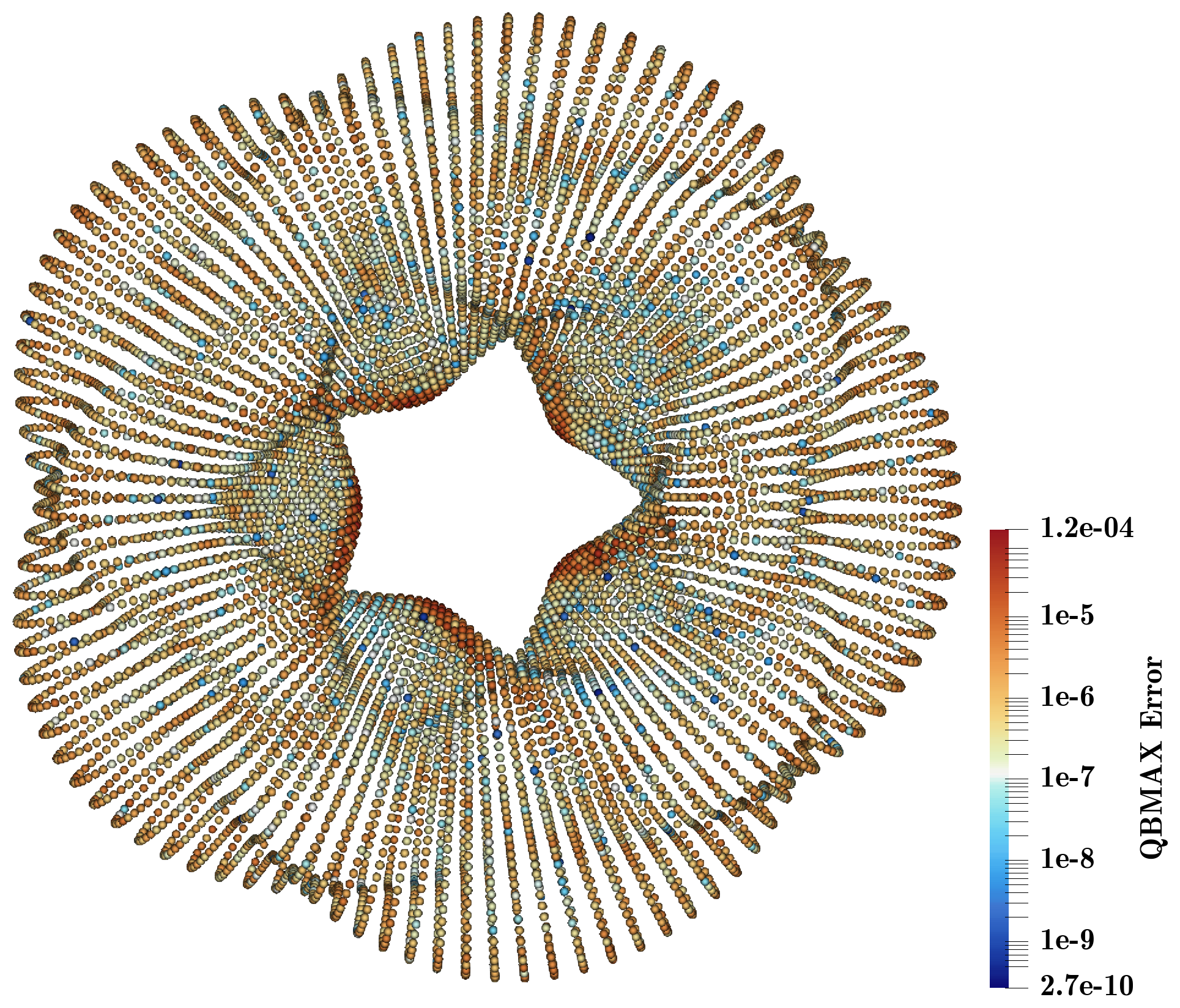}
        \caption{QBMAX}
    \end{subfigure}\hfill
    \begin{subfigure}[t]{0.48\textwidth}
        \centering \includegraphics[width=\textwidth]{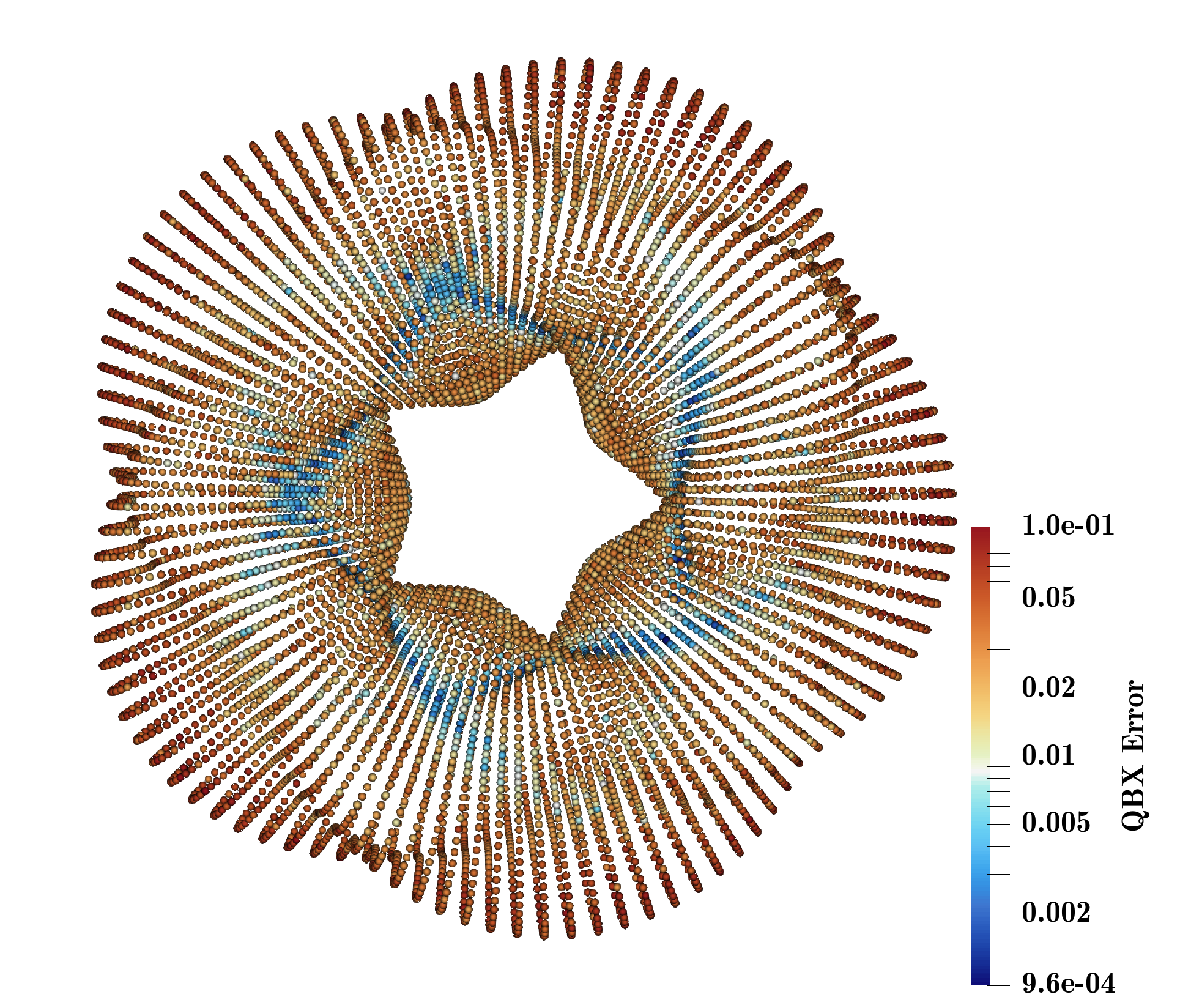}
        \caption{QBX}
    \end{subfigure}
    \caption{Error distribution of the jump relation for $\mathcal{D}_{k}\sigma$
    on the cruller surface for $k=40$.} \label{fig:error_distribution_cruller}
\end{figure}

 As shown in \Cref{fig:cruller_D_convergence}, convergence behavior for the jump
 relation of $\mathcal{D}_{k}\sigma$ can vary substantially between the two methods.
 \Cref{tab:cruller_jump_relations} summarizes the convergence rates for each
 method. While QBMAX achieves global convergence rates ranging from $3.5$ to
 $5.9$, QBX stagnates near $1.6$ when $k=80$, suggesting deterioration as $k$
 increases.
\begin{table}
    \centering
    \begin{subtable}[t]{0.9\textwidth}
       \centering \caption{$k = 10$}
       {\small\begin{tabular}[t]{|c|cc|cc|}
\hline
 $k = 10$ & \multicolumn{2}{c|}{QBMAX$\scriptstyle{(\mathcal{D}_k \sigma)}$} & \multicolumn{2}{c|}{QBX$\scriptstyle{(\mathcal{D}_k \sigma)}$} \\
\hline
$h$ & Error & EOC & Error & EOC \\
\hline
3.5e-01 & 4.2e-03 &   & 4.0e-02 &   \\
2.7e-01 & 1.6e-03 & 3.3 & 1.4e-02 & 3.8 \\
2.1e-01 & 7.2e-04 & 3.6 & 5.6e-03 & 4.0 \\
1.8e-01 & 3.6e-04 & 3.8 & 2.6e-03 & 4.2 \\
\hline
 & & 3.5 & & 4.0 \\
\hline
\end{tabular}
\begin{tabular}[t]{|c|}
\hline
\multicolumn{1}{|c|}{$\text{QBX} / \text{QBMAX}$} \\
Error Ratio \\
\hline
9.6 \\
8.5 \\
7.8 \\
7.2 \\
\hline
\end{tabular}}
    \end{subtable}

    \vspace{1em}

    \begin{subtable}[t]{0.9\textwidth}
       \centering \caption{$k = 20$}
       {\small\begin{tabular}[t]{|c|cc|cc|}
\hline
 $k = 20$ & \multicolumn{2}{c|}{QBMAX$\scriptstyle{(\mathcal{D}_k \sigma)}$} & \multicolumn{2}{c|}{QBX$\scriptstyle{(\mathcal{D}_k \sigma)}$} \\
\hline
$h$ & Error & EOC & Error & EOC \\
\hline
3.5e-01 & 2.7e-03 &   & 1.5e-01 &   \\
2.7e-01 & 1.0e-03 & 3.4 & 5.3e-02 & 3.6 \\
2.1e-01 & 4.4e-04 & 3.7 & 2.2e-02 & 3.9 \\
1.8e-01 & 2.2e-04 & 3.8 & 1.1e-02 & 4.1 \\
\hline
 & & 3.6 & & 3.8 \\
\hline
\end{tabular}
\begin{tabular}[t]{|c|}
\hline
\multicolumn{1}{|c|}{$\text{QBX} / \text{QBMAX}$} \\
Error Ratio \\
\hline
54.0 \\
52.2 \\
50.4 \\
48.2 \\
\hline
\end{tabular}}
    \end{subtable}

    \vspace{1em}

    \begin{subtable}[t]{0.9\textwidth}
       \centering \caption{$k = 40$}
       {\small\begin{tabular}[t]{|c|cc|cc|}
\hline
 $k = 40$ & \multicolumn{2}{c|}{QBMAX$\scriptstyle{(\mathcal{D}_k \sigma)}$} & \multicolumn{2}{c|}{QBX$\scriptstyle{(\mathcal{D}_k \sigma)}$} \\
\hline
$h$ & Error & EOC & Error & EOC \\
\hline
3.5e-01 & 2.1e-03 &   & 8.5e-01 &   \\
2.7e-01 & 4.9e-04 & 5.0 & 3.9e-01 & 2.7 \\
2.1e-01 & 2.3e-04 & 3.5 & 2.0e-01 & 3.1 \\
1.8e-01 & 1.2e-04 & 3.6 & 1.0e-01 & 3.4 \\
\hline
 & & 4.1 & & 3.0 \\
\hline
\end{tabular}
\begin{tabular}[t]{|c|}
\hline
\multicolumn{1}{|c|}{$\text{QBX} / \text{QBMAX}$} \\
Error Ratio \\
\hline
4.0e+02 \\
8.0e+02 \\
8.7e+02 \\
8.9e+02 \\
\hline
\end{tabular}}
    \end{subtable}

    \vspace{1em}

    \begin{subtable}[t]{0.9\textwidth}
       \centering \caption{$k = 80$}
       {\small\begin{tabular}[t]{|c|cc|cc|}
\hline
 $k = 80$ & \multicolumn{2}{c|}{QBMAX$\scriptstyle{(\mathcal{D}_k \sigma)}$} & \multicolumn{2}{c|}{QBX$\scriptstyle{(\mathcal{D}_k \sigma)}$} \\
\hline
$h$ & Error & EOC & Error & EOC \\
\hline
3.5e-01 & 4.7e-03 &   & 2.4e+00 &   \\
2.7e-01 & 7.8e-04 & 6.2 & 1.8e+00 & 1.1 \\
2.1e-01 & 2.1e-04 & 5.8 & 1.2e+00 & 1.7 \\
1.8e-01 & 7.8e-05 & 5.6 & 8.1e-01 & 2.2 \\
\hline
 & & 5.9 & & 1.6 \\
\hline
\end{tabular}
\begin{tabular}[t]{|c|}
\hline
\multicolumn{1}{|c|}{$\text{QBX} / \text{QBMAX}$} \\
Error Ratio \\
\hline
5.2e+02 \\
2.3e+03 \\
5.6e+03 \\
1.0e+04 \\
\hline
\end{tabular}}
    \end{subtable}
    \caption{Estimated order of convergence for jump relations of
    $\mathcal{D}_{k}\sigma$ on the cruller surface at $k = 10, 20, 40, 80$.}
    \label{tab:cruller_jump_relations}
\end{table}
Error patterns across the cruller surface in
\Cref{fig:error_distribution_cruller} reveal notable geometric dependencies. The
mean error for QBMAX is around $10^{-6}$, with the highest error occurring in
the inner region of the cruller surface, where high curvature is present. In
contrast, the QBX method exhibits a mean error around $10^{-2}$, with the
maximum error located in the outer region of the surface. These
curvature-dependent error distributions highlight the potential benefits of
local curvature correction, which will be explored in future work.
\section{Conclusions}
In this work, we have introduced QBMAX, a high-order method for the close
evaluation of layer potentials with rapidly decaying kernels, particularly
Yukawa potentials in 2D and 3D. By factoring out the leading exponential decay,
QBMAX reduces relative truncation error in the infinite-flat-boundary model and
generally reduces the degradation observed for QBX as $k$ increases in our
numerical tests. The magnitude of the improvement depends on the density,
geometry, and discretization parameters. Several directions for
future research remain open. First, extending the rigorous asymptotic analysis
beyond flat panels to general curved surfaces would provide further theoretical
foundations for the observed numerical performance. Second, integration with
fast algorithms such as the Fast Multipole Method
\cite{greengardNewVersionFast2002} would enable efficient large-scale
computations. Finally, the underlying principles of matched asymptotic
expansions can be naturally extended to other parameter-dependent PDEs,
including the modified biharmonic equation and modified Stokes flows. Results
along these directions will be presented in subsequent papers.

\section*{Acknowledgments}

The authors' research was supported by the National Science Foundation
under grants DMS-1654756, SHF-1911019, OAC-1931577, and DMS-2410943, and by the
US Department of Energy under award number DE-NA0002374,
as well as the Siebel School of
Computing and Data Science at the University of Illinois at Urbana-Champaign.
The work was further supported by Anthropic PBC under the `Claude for Scientists'
program.

A number of generative AI tools were used to assist in editing and reviewing this manuscript.
In addition, generative AI tools assisted with reviewing the experimental code,
verifying numerical results, and generating plots from computed data.
Every change made (to the manuscript and/or the experimental code) by an AI assistant
was reviewed by the authors before acceptance.
Models/services used include Qwen3.8-27B, OpenAI GPT-5.6 Terra and Sol,
and Anthropic Claude.

\appendix
\section{Integral identities}
\renewcommand{\thetheorem}{\Alph{section}.\arabic{theorem}}
\setcounter{theorem}{0}
In this section, we derive some closed-form expressions for integrals involving
the modified Bessel functions of the second kind.

\begin{lemma}\label{lem:flat-panel}
$\int_{\R} K_0(k\sqrt{s^{2} + r^{2}}) e^{i \omega s}ds =
\frac{\pi}{\sqrt{k^{2} + \omega^{2}}} \exp(-r \sqrt{k^{2} + \omega^{2}} )$ for
$k>0, \omega \in \mathbb{R} \text{ and } r \geq 0$.
\end{lemma}
\begin{myproof}
$K_0(k\sqrt{s^{2} + r^{2}})$ has an integral representation
\cite[\href{https://dlmf.nist.gov/10.32.E10}{(10.32.10)}]{DLMFNISTDigital}:
\[
K_0(k\sqrt{s^{2} + r^{2}}) = \frac{1}{2} \int_{0}^{\infty} \exp(-t - \frac{k^{2}(s^{2} + r^{2})}{4t}) \frac{dt}{t}.
\]
Then,
	\begin{align*}
		\int_{\R} K_0(k \sqrt{s^{2} + r^{2}}) e^{i \omega s} ds &=\frac{1}{2} \int_{\R} \int_{0}^{\infty} \exp(-t - \frac{k^{2}(s^{2} + r^{2})}{4t}) \frac{dt}{t} e^{i \omega s} ds\\
		&= \frac{1}{2} \int_{0}^{\infty} \exp(-t - \frac{k^{2} r^{2}}{4t}) \int_{\R} \exp(i \omega s - \frac{k^{2} s^{2}}{4t}) ds \frac{dt}{t}\\
		&= \frac{1}{2} \int_{0}^{\infty} \exp(-t - \frac{k^{2} r^{2}}{4t}) \sqrt{\frac{\pi}{k^{2} / 4t}}\exp(- \frac{\omega^{2}}{k^{2} / t}) \frac{dt}{t}\\
		&= \frac{\sqrt{\pi}}{k}\int_{0}^{\infty} \exp(-t(1+\frac{\omega^{2}}{k^{2}}) - \frac{k^{2} r^{2}}{4t}) \frac{dt}{\sqrt{t}}\\
		&= \frac{\sqrt{\pi}}{k} \int_{0}^{\infty} \exp(-u - \frac{k^{2}r^{2}}{4}\frac{(1+ \frac{\omega^{2}}{k^{2}})}{u}) (1+ \frac{\omega^{2}}{k^{2}})^{-1} \sqrt{1+ \frac{\omega^{2}}{k^{2}}} \frac{du}{\sqrt{u}}\\
		&= \frac{\sqrt{\pi}}{k} (1+ \frac{\omega^{2}}{k^{2}})^{-\frac{1}{2}} \int_{0}^{\infty}\exp(-u - \frac{k^{2}r^{2} + r^{2} \omega^{2}}{4u}) \frac{du}{\sqrt{u}} \\
		&= \sqrt{\pi} (k^{2} + \omega^{2})^{-\frac{1}{2}} [2 (\frac{1}{2}|r| \sqrt{k^{2} + \omega^{2}})^{\frac{1}{2}} K_{-\frac{1}{2}}(|r| \sqrt{k^{2} + \omega^{2}}) ] \\
		&= \sqrt{\pi} (k^{2} + \omega^{2})^{-\frac{1}{2}} [2 (\frac{1}{2} |r| \sqrt{k^{2} + \omega^{2}})^{\frac{1}{2}} \sqrt{\frac{\pi}{2}} \frac{\exp(-|r| \sqrt{k^{2} + \omega^{2}})}{\sqrt{|r| \sqrt{k^{2} + \omega^{2}}}}] \\
		&= \frac{\pi}{\sqrt{k^{2} + \omega^{2}}} \exp(-|r| \sqrt{k^{2} + \omega^{2}}),
	\end{align*}
	where the second-to-last equality follows from
	$K_{-1/2}(z)=\sqrt{\pi/(2z)}e^{-z}$
	\cite[\href{https://dlmf.nist.gov/10.39.E2}{(10.39.2)}]{DLMFNISTDigital}.
\end{myproof}

\begin{lemma}\label{lem:disk}
Suppose $\sigma_{m}(\theta) = e^{i m \theta}$ is a layer density on $\mathbb{T}$
for some $m \in \Z$ and $k>0$.
\[
	\mathcal{S}_{k} \sigma(\b x) = \frac{1}{2\pi}\int_{\mathbb{T}} K_{0}(k|\b x-\b y(\theta)|)  \sigma_m(\theta) d\theta = \begin{cases} e^{i m \theta_{\b x} } I_m (k|\b x|) K_m(k), & |\b x| \leq 1, \\
		e^{i m \theta_{\b x}} I_m(k) K_m(k|\b x|), & |\b x| > 1, \end{cases}
\]
where $\b y(\theta) = (\cos(\theta), \sin (\theta))$.
\end{lemma}

\begin{myproof}
Suppose $\b x = r_{\b x} (\cos \theta_{\b x},\sin \theta_{\b x})$ and $\b y = (\cos \theta, \sin
\theta)$ for some $1 > r_{\b x} \geq 0$. By Graf's addition theorem for
modified Bessel functions
\cite[\href{https://dlmf.nist.gov/10.44.ii}{10.44(ii)}]{DLMFNISTDigital},
\[
K_0(k|\b x-\b y|) = \sum_{l \in \Z} K_{l} (k|\b y|) I_{l}(k|\b x|) e^{il (\theta -\theta_{\b x})},
\qquad |\b x| < |\b y|.
\]
Hence,
\[
\mathcal{S}_{k} \sigma(\b x) =  \frac{1}{2\pi} \int_{0}^{2\pi} \sum_{l \in \Z} K_{l}(k) I_{l}(k|\b x|) e^{il (\theta - \theta_{\b x})} e^{i m \theta } d \theta = e^{i m \theta_{\b x}} I_m(k|\b x|) K_m(k).
\]
Similarly, for $|\b x| >1$,
\[
\mathcal{S}_{k} \sigma(\b x) = e^{im \theta_{\b x}} I_m(k) K_m(k|\b x|).
\]
The continuity of $\mathcal{S}_{k}\sigma$ on $\mathbb{R}^2$ implies
$\mathcal{S}_{k}\sigma(\b x)=e^{im \theta_{\b x}}I_m(k)K_m(k)$ for $|\b x| =1$.

\end{myproof}

\begin{lemma}\label{lem:dis_D}
	Suppose $\sigma_m (\theta)= e^{im \theta}$ is a layer density on $\mathbb{T}$
	for some $m \in \Z$ and $k>0$, then
\[
	\mathcal{D}_{k} \sigma(\b x) = \frac{1}{2\pi}\int_{\mathbb{T}} \frac{\partial }{\partial \b \nu_{\b y}}K_{0}(k|\b x-\b y(\theta)|) \sigma_{m}(\theta) d\theta = \begin{cases} k e^{i m \theta_{\b x}} I_{m}(k|\b x|) K'_{m}(k), & |\b x| < 1, \\
			k e^{ i m \theta_{\b x}}K_{m}(k|\b x|) I_m'(k), & |\b x| > 1, \end{cases}
\]
where $\b y(\theta) = (\cos(\theta), \sin(\theta))$.
\end{lemma}
\begin{myproof}
Suppose $\b x = r_{\b x}(\cos \theta_{\b x}, \sin \theta_{\b x})$ and $\b y=r (\cos \theta,
\sin \theta)$ for some $r> r_{\b x} \geq 0$. By Graf's addition theorem for
modified Bessel functions
\cite[\href{https://dlmf.nist.gov/10.44.ii}{10.44(ii)}]{DLMFNISTDigital}, $K_0(k|\b x-\b y|)$ can be
represented as
\[
K_0(k|\b x-\b y|) = \sum_{l \in \Z} K_l (k|\b y|) I_{l}(k|\b x|) e^{il(\theta -\theta_{\b x})},
\qquad |\b x| < |\b y|.
\]
\[
\begin{split}
\langle \nabla_{\b y} [K_{l}(k|\b y|) e^{il \theta}], \frac{\b y}{|\b y|} \rangle &= k e^{il \theta} K'_{l}(k|\b y|) + K_{l}(k|\b y|)\langle \nabla_{\b y} (\frac{y_1 + i y_2}{|\b y|})^{l}, \frac{\b y}{|\b y|}  \rangle \\
& = k e^{il \theta} K_{l}'(k|\b y|).
\end{split}
\]
The second inner product vanishes because $e^{il\theta}$ is constant along
radial rays, so its directional derivative in the direction
$\b y/|\b y|$ is zero. Hence,
\[
\mathcal{D}_{k} \sigma(\b x) = \frac{1}{2\pi}\int _{0}^{2\pi} \sum_{l \in \Z} kK'_{l}(k) I_{l}(k|\b x|) e^{il (\theta - \theta_{\b x})} e^{im \theta} d\theta = k e^{i m \theta_{\b x}}K_{m}'(k) I_{m}(k|\b x|).
\]
The case when $|\b x| > 1$ can be proved similarly with $\mathcal{D}_{k}\sigma (\b x)
= k e^{i m \theta_{\b x}} K_{m}(k|\b x|) I_m'(k)$.

\end{myproof}

\section{Derivative bounds for \texorpdfstring{$K_0$}{K0}}\label[appendix]{sec:appendix-derivative-bounds}

Recall the series expansion of the kernel $K_0(z)$ \cite[\href{https://dlmf.nist.gov/10.31.E2}{(10.31.2)}]{DLMFNISTDigital}:
\begin{equation*}
	\begin{split}
	K_0(z) = \ln\!\bigl(\tfrac{2}{z}\bigr)
		&+ \ln\!\bigl(\tfrac{2}{z}\bigr)\sum_{n=1}^{\infty}\frac{1}{(n!)^2}\Bigl(\frac{z^2}{4}\Bigr)^{n}
		- \gamma\Bigl(1 + \sum_{n=1}^{\infty}\frac{1}{(n!)^2}\Bigl(\frac{z^2}{4}\Bigr)^{n}\Bigr) \\
		&+ \sum_{n=1}^{\infty}\!\Bigl(1+\tfrac12+\cdots+\tfrac1n\Bigr)\frac{1}{(n!)^2}\Bigl(\frac{z^2}{4}\Bigr)^{n}.
	\end{split}
\end{equation*}
For $m \in \mathbb{N}_0$, the singularity of $\frac{\partial ^{m}}{\partial z^{m}} K_0(z)$ is of order $\mathcal{O}(\max\{z^{-m}, |\log z|\})$. Hence, for $0 < z < 1$,
$
   |K_{0}^{(m)}(z)| \leq C_m \max\{z^{-m}, |\log z|\}
$
for some constant $C_m$ depending on the derivative order $m$.

Using the integral representation of $K_0(z)$ \cite[\href{https://dlmf.nist.gov/10.32.E9}{(10.32.9)}]{DLMFNISTDigital} for $m \in \mathbb{N}_0$,
\[
K_{0}^{(m)}(z) = \int_{0}^{\infty }(-\cosh(t))^m e^{-z \cosh(t)} dt, \qquad z > 0.
\]
Then, for $z \geq 1$,
\begin{equation*}
\begin{split}
|K_{0}^{(m)}(z)| &\leq  e^{-z} \int_0 ^{\infty} (\cosh(t))^{m} e^{-z (\cosh(t) - 1)} dt\\
&= e^{-z} \int_{0}^{\infty} \left[(\cosh(t))^{m} e^{-\frac{z}{2} (\cosh(t) - 1)}\right] e^{-\frac{z}{2} (\cosh(t) - 1)} dt \\
& \leq e^{-z} \int_{0}^{\infty} \left[(\cosh(t))^{m} e^{-\frac{1}{2} (\cosh(t) - 1)}\right] e^{-\frac{z}{2} (\cosh(t) - 1)} dt \\
&\leq C_m e^{-z}  \int_{0}^{\infty} e^{-\frac{z}{2} (\cosh(t) - 1)} dt,
\end{split}
\end{equation*}
for some constant $C_m$ depending on the derivative order $m$.
Using the series expansion for $\cosh(t) - 1$,
\[
    \cosh(t) - 1 = \sum_{n=1}^{\infty} \frac{t^{2n}}{(2n)!},
\]
we conclude that
\[
   |K_0^{(m)}(z)| \leq   C_m e^{-z} \int_{0}^{\infty} e^{-\frac{z}{4} t^{2}} dt = C_m e^{-z} \sqrt{\frac{\pi}{z}}, \qquad z \geq 1.
\]
Combining the above estimates, for $m \in \mathbb{N}_0$,
\[
    |K_0^{(m)}(z)| \leq C_m \begin{cases}
        \max\{z^{-m}, |\log z|\}, &  0 < z < 1, \\
        z^{-\frac{1}{2}} e^{-z}, & z \geq 1,
    \end{cases}
\]
for some $C_m > 0$ depending only on the derivative order $m$.

\section{Quadrature error analysis for the QBX expansion}\label[appendix]{sec:appendix-quad-error}
In this appendix, and in support of the argument in \Cref{sec:quad-error}, we estimate the quadrature error of the QBX (not QBMAX) line expansion of the single-layer potential
\[
    \mathcal{S}_k\sigma(\b x) = \int_{\partial\Omega} \mathcal{K}(\b x,\b y)\,\sigma(\b y)\,dS_{\b y}, \qquad \b x \in \partial\Omega,
\]
with kernel $\mathcal{K}(\b x,\b y) = \frac{1}{2\pi}K_0(k|\b x-\b y|)$ for on-surface evaluation. We assume that $\sigma \in C^{\infty}(\partial\Omega)$ and that the $C^{\infty}$-smooth boundary $\partial\Omega$ is discretized into equal-size, non-overlapping panels $\Gamma_i$ of arc length $h$, parameterized by arc length as $\b \gamma_i: [0,h] \to \Gamma_i$. We further assume that these panel parameterizations are restrictions of a fixed smooth global parametrization, so their derivatives are bounded uniformly under refinement.
We further make liberal use of small-argument approximations of the kernel (cf.~\Cref{sec:appendix-derivative-bounds}).
The derivation parallels that of~\cite{epsteinConvergenceLocalExpansions2013a}, where an analogous result is shown for the Laplace kernel.

Fix a target $\b x \in \partial\Omega$ and an expansion center $\b c_{\b x}$ at distance $r = |\b x - \b c_{\b x}|$, and assume the standard QBX clearance condition
$\overline{B_r(\b c_{\b x})}\cap\partial\Omega=\{\b x\}$. In particular,
$|\b c_{\b x}-\b y|\geq r$ for every $\b y\in\partial\Omega$. The $p^{\text{th}}$-order QBX line expansion of $\mathcal{S}_k\sigma$ at $\b x$ is defined as
\[
    \mathcal{T}_p[\mathcal{S}_k\sigma\circ \b l](1) = \sum_{m=0}^p \frac{1}{m!}\int_{\partial\Omega} \partial_\tau^m \mathcal{K}(\b l(\tau), \b y)\big\rvert_{\tau=0} \sigma(\b y) dS_{\b y},
\]
where $\b l(\tau)= \b c_{\b x} + \tau(\b x - \b c_{\b x})$. Following \cite[Theorem~1]{klocknerQuadratureExpansionNew2013} and \cite{epsteinConvergenceLocalExpansions2013a}, we estimate the quadrature error.

Define $D_m(\b y) := \partial_\tau^m \mathcal{K}(\b l(\tau), \b y)\big\rvert_{\tau=0}$ and
$F_{i,m}(s):= D_m(\b \gamma_i(s)) \sigma(\b \gamma_i(s))$ for $s \in [0,h]$.
Let $Q_{h}$ be a $q$-point Gauss--Legendre quadrature rule on $[0,h]$. Then, by the quadrature error estimate
\cite{davisMethodsNumericalIntegration2014} and the Leibniz rule,
\begin{align*}
   \left|\int_{0}^{h} F_{i,m}(s) ds - Q_{h}(F_{i,m})\right|
   &\leq \frac{(q!)^{4}}{(2q + 1)[(2q)!]^{3}} h^{2q+1}\sup_{s \in [0,h]}|F_{i,m}^{(2q)}(s)| \\
   &= \frac{(q!)^{4}}{(2q + 1)[(2q)!]^{3}} h^{2q+1}\sup_{s \in [0,h]} \left| \sum_{j=0}^{2q} \binom{2q}{j} \partial_{s}^{j} D_m(\b \gamma_i(s))\,\partial_{s}^{2q-j} \sigma(\b \gamma_i(s)) \right|.
\end{align*}
Since $\b l(\tau) = (\b x - \b c_{\b x})\tau + \b c_{\b x}$,
using multi-index notation $\b{\alpha} = (\b{\alpha}_1, \b{\alpha}_2) \in \mathbb{N}_0^2$ with $|\b{\alpha}| = \b{\alpha}_1 + \b{\alpha}_2$, $\b{\alpha}! = \b{\alpha}_1!\,\b{\alpha}_2!$, and $\b z^{\b{\alpha}} = z_1^{\b{\alpha}_1} z_2^{\b{\alpha}_2}$, $D_m(\b y)$ can be simplified as
\[
    D_m(\b y) = \left. ( (\b x - \b c_{\b x}) \cdot \nabla_{\b z})^m \mathcal{K}(\b z, \b y) \right\rvert_{\b z = \b c_{\b x}}
    = \sum_{\b{\alpha} \in \mathbb{N}_0^2,\,|\b{\alpha}| = m} \frac{m!}{\b{\alpha}!} (\b x - \b c_{\b x})^{\b{\alpha}} \partial^{\b{\alpha}}_{\b z} \mathcal{K}(\b z, \b y)\big\rvert_{\b z = \b c_{\b x}}.
\]
Using Fa\`a di Bruno's formula \cite{constantineMultivariateFAADi1996}, $\partial_{s}^{j} D_m(\b \gamma_i(s))$ can be expressed as
\begin{equation*}
\begin{split}
   \partial_{s}^{j} D_m(\b \gamma_i(s)) &= \sum_{\b{\beta} \in \mathbb{N}_0^2,\,|\b{\beta}| \leq j} P_{j,\b{\beta}} \partial_{\b y}^{\b{\beta}} D_m(\b y)\big\rvert_{\b y = \b \gamma_i(s)}\\
   &= \sum_{\b{\beta} \in \mathbb{N}_0^2,\,|\b{\beta}| \leq j} P_{j,\b{\beta}} \sum_{\b{\alpha} \in \mathbb{N}_0^2,\,|\b{\alpha}| = m} \frac{m!}{\b{\alpha}!} (\b x - \b c_{\b x})^{\b{\alpha}} \partial_{\b y}^{\b{\beta}} \partial^{\b{\alpha}}_{\b z} \mathcal{K}(\b z, \b y)\big\rvert_{\b z = \b c_{\b x}, \b y = \b \gamma_i(s)},
\end{split}
\end{equation*}
where $P_{j,\b{\beta}}(\b \gamma_i'(s), \b \gamma_i''(s), \dots, \b
\gamma_i^{(j)}(s))$ is a polynomial in the derivatives of $\b \gamma_i$.
For $m+|\b{\beta}|>0$, repeated application of the chain rule to the radial
function $K_0(k|\b z-\b y|)$, together with the derivative estimates in
\Cref{sec:appendix-derivative-bounds}, gives
\[
 |\partial_{\b y}^{\b{\beta}}\partial_{\b z}^{\b{\alpha}}
 \mathcal K(\b z,\b y)|
 \leq C(k,m,|\b{\beta}|)|\b z-\b y|^{-m-|\b{\beta}|},
 \qquad |\b{\alpha}|=m,
\]
for $0<|\b z-\b y|\leq1$. Indeed, each Cartesian derivative either
differentiates a radial derivative of $K_0$ or a factor
$(z_j-y_j)/|\b z-\b y|$, and hence contributes at most one additional inverse
power of $|\b z-\b y|$. For larger separation the exponential-decay estimate
in \Cref{sec:appendix-derivative-bounds} is bounded and may be absorbed into
the same constant. Consequently,
\[
\left|(\b x-\b c_{\b x})^{\b{\alpha}}
\partial_{\b y}^{\b{\beta}}\partial_{\b z}^{\b{\alpha}}
\mathcal{K}(\b c_{\b x},\b y)\right|
\leq C(k, m, |\b{\beta}|)r^m
|\b c_{\b x}-\b y|^{-m-|\b{\beta}|}
\leq C(k, m, |\b{\beta}|)r^{-|\b{\beta}|}.
\]
Since $q$ is a fixed positive integer and
$|\b{\beta}|\leq j\leq 2q$,
$r^{-|\b{\beta}|}=\mathcal O(r^{-2q})$ as $r\to0$.
For $m=|\b{\beta}|=0$,
\Cref{sec:appendix-derivative-bounds} gives
\[
|\mathcal{K}(\b c_{\b x},\b y)|
\leq C(k)|\log r|=\mathcal O(r^{-2q})
\quad\text{as }r\to0.
\]
Thus, for sufficiently small $r>0$,
\[
\sup_{s \in [0,h]} |\partial_s^j D_m(\b \gamma_i(s))|
\leq C(k, m, q, \partial\Omega) \frac{m!}{r^{2q}},
\qquad j=0,\ldots,2q,
\]
for some constant $C(k, m, q, \partial\Omega)$. Therefore, the quadrature
error for the $m$-th term in the QBX expansion can be bounded by
\[
   \frac{1}{m!}| \int_{0}^{h} F_{i,m}(s) ds - Q_{h}(F_{i,m})| \leq C(k, m, q, \partial\Omega) \frac{(q!)^{4}}{(2q + 1)[(2q)!]^{3}} h^{2q+1} \frac{2^{2q}}{r^{2q}} \left\|\sigma\right\|_{C^{2q}(\partial\Omega)}.
\]
By Stirling's approximation $\sqrt{2\pi}n^{n+1/2}e^{-n} \leq n! \leq 2\sqrt{\pi}n^{n+1/2}e^{-n}$:
\[
    (q!)^4 \leq (2\sqrt{\pi})^4 q^{4q+2} e^{-4q} \text{ and } 
    [(2q)!]^3 \geq (\sqrt{2\pi})^3 (2q)^{6q+3/2} e^{-6q}.
\]
Hence,
\[
   \frac{(q!)^4 2^{2q}}{(2q+1)[(2q)!]^3} \leq \frac{2\sqrt{\pi q}}{2q+1} \left(\frac{e}{4q}\right)^{2q}.
\]
Therefore,
\[
      \sum_{m=0}^{p} \frac{1}{m!}\big| \int_{0}^{h} F_{i,m}(s) ds - Q_{h}(F_{i,m})\big| \leq C(k, p, q, \partial\Omega) \frac{h^{2q+1}}{4^{2q} r^{2q}} \|\sigma\|_{C^{2q}(\partial\Omega)}.
\]
Summing over all panels $\Gamma_i$, the quadrature error for the $p^{\text{th}}$ order QBX expansion can be bounded by
\[
   \left| \sum_{i} \sum_{m=0}^{p} \frac{1}{m!} \left[\int_{0}^{h} F_{i,m}(s) ds - Q_{h}(F_{i,m})\right]\right| 
   \leq C(k, p, q, \partial\Omega) \left(\frac{h}{4r}\right)^{2q} \|\sigma\|_{C^{2q}(\partial\Omega)},
\]
for some constant $C(k, p, q, \partial\Omega)$ depending on $k$, $p$, $q$ and $\partial\Omega$.

\bibliographystyle{elsarticle-num}
\bibliography{reference}

\end{document}